\documentclass[11pt]{article}
\usepackage{geometry}                
\usepackage{amsmath,amssymb}
\usepackage{amsthm}
\usepackage{hyperref}
\usepackage{graphicx,epsfig,color}
\usepackage{booktabs,bm,multirow}
\usepackage{cases}

\newtheorem{theorem}{Theorem}

\newtheorem{remark}[theorem]{Remark}
\newtheorem{lemma}[theorem]{Lemma}

\input{siam.sty}

\title{Randomized extended block Kaczmarz for solving least squares}
\author{Kui Du\thanks{School of Mathematical Sciences, Xiamen University, Xiamen 361005, China ({\tt kuidu@xmu.edu.cn}).},\quad Wu-Tao Si\thanks{School of Mathematical Sciences, Xiamen University, Xiamen 361005, China ({\tt 19020171152496@stu.xmu.edu.cn}).},\quad  Xiao-Hui Sun\thanks{School of Mathematical Sciences, Xiamen University, Xiamen 361005, China ({\tt 19020190154621@stu.xmu.edu.cn}).}} 
\date{}                                           

\begin{document}
\maketitle

\begin{abstract} Randomized iterative algorithms have recently been proposed to solve large-scale linear systems. In this paper, we present a simple randomized extended block Kaczmarz algorithm that exponentially converges in the mean square to the unique minimum $\ell_2$-norm least squares solution of a given linear system of equations. The proposed algorithm is pseudoinverse-free and therefore different from the projection-based randomized double block Kaczmarz algorithm of Needell, Zhao, and Zouzias. We emphasize that our method works for all types of linear systems (consistent or inconsistent, overdetermined or underdetermined,  full-rank or rank-deficient).  Moreover, our approach can utilize efficient implementations on distributed computing units, yielding remarkable improvements in computational time. Numerical examples are given to show the efficiency of the new algorithm.  
\vspace{.5mm} 


\end{abstract}

\section{Introduction} 

The Kaczmarz method \cite{kaczmarz1937angen} is a simple iterative method for solving a linear systems of equations $${\bf Ax=b},\quad \mbf A\in\mbbr^{m\times n}, \quad \mbf b\in\mbbr^m.$$ Due to its simplicity and numerical performance, the Kaczmarz method has found many applications in many fields, such as computer tomography \cite{natterer2001mathe,kak2001princ,herman2009funda}, image reconstruction \cite{popa2004kaczm,herman2008image}, digital signal processing \cite{byrne2004unifi,lorenz2015spars}, etc.
At each step, the method projects the current iterate onto one hyperplane defined by a row of the system.  More precisely, assuming that the $i$th row $\mbf A_{i,:}$ has been selected at the $k$th iteration, then the $k$th estimate vector $\mbf x^k$ is obtained by $$\mbf x^k=\mbf x^{k-1}-\alpha_k\frac{\mbf A_{i,:}\mbf x^{k-1}-\mbf b_i}{\mbf A_{i,:}(\mbf A_{i,:})^\rmt}(\mbf A_{i,:})^\rmt,$$ where $(\mbf A_{i,:})^\rmt$ denotes the transpose of $\mbf A_{i,:}$, $\mbf b_i$ is the $i$th component of $\mbf b$, and $\alpha_k$ is a stepsize. Numerical experiments show that using the rows of the coefficient matrix in the Kaczmarz method in random order, rather than in their given order, can often greatly improve the convergence \cite{herman1993algeb,natterer2001mathe}. In a seminal paper \cite{strohmer2009rando}, Strohmer and Vershynin proposed a randomized  Kaczmarz (RK) algorithm  which exponentially converges in expectation to the solutions of consistent, overdetermined, full-rank linear systems. The convergence result was extended and refined in {various directions including inconsistent \cite{leventhal2010rando,needell2010rando,zouzias2013rando,dumitrescu2015relat,petra2016singl,needell2016stoch,haddock2019motzk}, underdetermined or rank-deficient linear systems \cite{ma2015conve,gower2015rando,razaviyayn2019linea,du2019doubl}, ridge regression problems \cite{hefny2017rows,liu2019varia}, linear feasibility problems \cite{loera2017sampl}, convex feasibility problems \cite{necoara2019rando}, block variants \cite{needell2014paved,needell2015rando,necoara2019faste}, acceleration strategies \cite{liu2016accel,richtarik2017stoch,bai2018greed,bai2018relax,bai2019greed,bai2019parti,zhang2019new}, and many others \cite{bai2018conve,popa2018conve,du2019tight,du2019new,haddock2019rando}.}

Let $\mbf A^\dag$ denote the Moore-Penrose pseudoinverse\footnote{Every $m\times n$ matrix $\mbf A$ has a unique Moore-Penrose pseudoinverse. In particular, in this paper we will use the following property of the pseudoinverse: $\mbf A^\rmt=\mbf A^\rmt\mbf A\mbf A^\dag$.}  \cite{ben2003gener} of $\mbf A$. In this paper, we are interested in the vector $\mbf A^\dag\mbf b$. Here we would like to make clear what $\mbf A^\dag \mbf b$ stands for different types of linear systems (see \cite{ben2003gener,golub2013matri}):
\bit
\item[(1)] If $\bf Ax=b$ is consistent with full-column rank $\mbf A$, i.e.,  $\rank(\mbf A)=n$, then $\mbf A^\dag\mbf b$ is the unique solution. In this case, we have $m\geq n$ and  the linear system is overdetermined when $m>n$. 
\item[(2)] If $\bf Ax=b$ is consistent with $\rank(\mbf A)<n$, then $\mbf A^\dag\mbf b$ is the unique minimum $\ell_2$-norm solution. In this case, we have $m\geq n$ or $m<n$, and  the linear system is overdetermined (resp. underdetermined) when $m>n$ (resp. $m<n$). The matrix $\mbf A$ can be of full-row  rank, i.e., $\rank(\mbf A)=m$, or rank-deficient, i.e., $\rank(\mbf A)<m$. 
\item[(3)] If $\bf Ax=b$ is inconsistent with $\rank(\mbf A)=n$, then $\mbf A^\dag\mbf b$ is the unique least squares solution. In this case, we have $m\geq n$ and  the linear system is overdetermined when $m>n$.
\item[(4)] If $\bf Ax=b$ is inconsistent with $\rank(\mbf A)<n$, then $\mbf A^\dag\mbf b$ is the unique  minimum $\ell_2$-norm  least squares solution. In this case, we have $m\geq n$ or $m<n$, and  the linear system is overdetermined (resp. underdetermined) when $m>n$ (resp. $m<n$). The matrix $\mbf A$ can be of full-row  rank, i.e., $\rank(\mbf A)=m$, or rank-deficient, i.e., $\rank(\mbf A)<m$. 
\eit

If ${\bf Ax = b}$ is inconsistent, Needell \cite{needell2010rando} showed that RK does not converge to $\bf A^\dag b$. To resolve this problem, Zouzias and Freris \cite{zouzias2013rando} proposed a randomized extended Kaczmarz (REK) algorithm, which uses RK twice \cite{liu2016accel,du2019tight} at each iteration and exponentially converges in the mean square to $\bf A^\dag b$. More precisely, assuming that the $j$th column $\mbf A_{:,j}$ and the $i$th row $\mbf A_{i,:}$ have been selected at the $k$th iteration, REK generates two vectors $\mbf z^{k}$ and $\mbf x^{k}$ via two RK updates (one for $\bf A^\rmt z= 0$ from $\mbf z^{k-1}$ and the other for ${\bf Ax = b-z}^k$ from $\mbf x^{k-1}$): \begin{align*}\mbf z^k &= \mbf z^{k-1}-\dsp\frac{(\mbf A_{:,j})^\rmt\mbf z^{k-1}}{(\mbf A_{:,j})^\rmt\mbf A_{:,j}}\mbf A_{:,j},\\ {\bf x}^k &= {\bf x}^{k-1}-\dsp\frac{\mbf A_{i,:}\mbf x^{k-1}-\mbf b_i+\mbf z^k_i}{{\bf A}_{i,:}({\bf A}_{i,:})^\rmt}(\mbf A_{i,:})^\rmt.\end{align*} For general linear systems (consistent or inconsistent, full-rank or rank-deficient), the vector $\mbf x^k$ generated by REK exponentially converges to $\mbf A^\dag\mbf b$ if $\mbf z^0\in\mbf b+\ran(\mbf A)$ and $\mbf x^0\in\ran({\mbf A^\rmt})$ \cite{liu2016accel,du2019tight}.
To accelerate the convergence, the following projection-based block variants \cite{needell2014paved,needell2015rando} of RK and REK were developed. For a subset $\mcali\subset\{1,2,\ldots,m\}$ and a subset $\mcalj\subset\{1,2,\ldots,n\}$, denote by $\mbf A_{\mcali,:}$ and $\mbf A_{:,\mcalj}$ the row submatrix of $\mbf A$ indexed by $\mcali$ and the column submatrix of $\mbf A$ indexed by $\mcalj$, respectively.
Assuming that the subset $\mcali_i$ has been selected at the $k$th iteration, the randomized block Kaczmarz (RBK) algorithm \cite{needell2014paved} generates the $k$th estimate $\mbf x^k$ via $$\mbf x^k=\mbf x^{k-1}-(\mbf A_{\mcali_i,:})^\dag(\mbf A_{\mcali_i,:}\mbf x^{k-1}-\mbf b_{\mcali_i}).$$  Assuming that the subsets $\mcalj_j$ and $\mcali_i$ have been selected at the $k$ iteration, the randomized double block Kaczmarz (RDBK) algorithm \cite{needell2015rando} generates the $k$th estimate $\mbf x^k$ via \begin{align*} \mbf z^k &= \mbf z^{k-1}-\mbf A_{:,\mcalj_j}(\mbf A_{:,\mcalj_j})^\dag\mbf z^{k-1},\\
{\bf x}^k &= {\bf x}^{k-1}-(\mbf A_{\mcali_i,:})^\dag(\mbf A_{\mcali_i,:}\mbf x^{k-1}-\mbf b_{\mcali_i}+\mbf z^k_{\mcali_i}).\end{align*}
Numerical experiments demonstrate that the convergence can be significantly accelerated if appropriate blocks of the coefficient matrix are used. {The main drawback of projection-based block methods is that  they are not adequate for distributed implementations}. 

Recently, Necoara \cite{necoara2019faste} proposed a randomized average block Kaczmarz (RABK) algorithm for consistent linear systems, which takes a convex combination of { several RK updates (i.e., the projections of the current iterate onto several hyperplanes)} as a new direction with some stepsize. Assuming that the subset $\mcali$ has been selected at the $k$th iteration, RABK generates the $k$th estimate $\mbf x^k$ via  \beq\label{rabk}\mbf x^k=\mbf x^{k-1}-\alpha_k\l(\sum_{i\in\mcali}\omega_i^k\frac{\mbf A_{i,:}\mbf x^{k-1}-\mbf b_i}{\mbf A_{i,:}(\mbf A_{i,:})^\rmt}(\mbf A_{i,:})^\rmt\r),\eeq where the weights $\omega_i^k\in[0,1]$ such that $\sum_{i\in\mcali}\omega_i^k=1$, and the stepsize $\alpha_k\in(0,2)$. The convergence analysis reveals that RABK is extremely effective when it is given a good sampling of the rows into well-conditioned blocks. { A block version of RABK (i.e., parallel randomized block Kaczmarz), which takes a convex combination of the RBK updates, was proposed and studied by Richt\'{a}rik and Tak\'{a}\v{c} \cite{richtarik2017stoch}.}  Shortly afterwards, Du and Sun \cite{du2019doubl} proposed a doubly stochastic block Gauss-Seidel (DSBGS) algorithm, which randomly chooses a submatrix of the coefficient matrix at each iteration. Assuming that the subsets $\mcali$ and $\mcalj$ have been selected at the $k$th iteration, DSBGS generates the $k$th estimate $\mbf x^k$ via $$\mbf x^k=\mbf x^{k-1}-\alpha_k \frac{\mbf I_{:,\mcalj}(\mbf A_{\mcali,\mcalj})^\rmt(\mbf I_{:,\mcali})^\rmt}{\|\mbf A_{\mcali,\mcalj}\|_\rmf^2} (\mbf A\mbf x^{k-1}-\mbf b),$$ where $\mbf I$ denotes the identity matrix, $\mbf A_{\mcali,\mcalj}$ denotes the submatrix that lies in the rows indexed by $\mcali$ and the columns indexed by $\mcalj$, and $\|\cdot\|_\rmf$ is the Frobenius norm. Exponential convergence of DSBGS for consistent linear systems was proved. By setting $\mcali\subset\{1,2,\ldots,m\}$ and $\mcalj=\{1,2,\ldots,n\}$, DSBGS recovers a special case of RABK, i.e., RABK with weight $$\omega_i^k=\frac{\mbf A_{i,:}(\mbf A_{i,:})^\rmt}{\|\mbf A_{\mcali,:}\|_\rmf^2},\qquad i\in\mcali.$$  Note that both RABK and DSBGS are { very easy to implement on distributed computing units}, yielding remarkable improvements in computational time. We emphasize that convergence results in the mean square of RABK and DSBGS are obtained only for consistent linear systems.

In this paper, { based on the REK algorithm and the RABK algorithm, we present a simple} randomized extended block Kaczmarz (REBK) algorithm that exponentially converges in the mean square to the unique minimum $\ell_2$-norm (least squares) solution of a given general linear system (full-rank or rank-deficient, overdetermined or underdetermined, consistent or inconsistent). Our method is different from those projection-based block methods, for example, those in  \cite{elfving1980block,arioli1992block,bramley1992row,popa1998exten,needell2014paved,needell2015rando,duff2015augme}. At each step, REBK, as a direct extension of REK, uses two {special RABK (which also can be viewed as special DSBGS)} updates (one for $\bf A^\rmt z= 0$ from $\mbf z^{k-1}$ and the other for ${\bf Ax = b-z}^k$ from $\mbf x^{k-1}$; see Section 2 for details). Compared with REK, REBK { usually has a better convergence rate and} can exploit the high-level basic linear algebra subroutine ({\tt BLAS2}), even fast matrix-vector multiplies (for example, if submatrices of $\mbf A$ have circulant or Toeplitz structures, then the Fast Fourier Transform technique can be used), and therefore could be more efficient. { Compared with RDBK, REBK can be implemented on distributed computing units}.  We refer the reader to \cite{needell2014paved,necoara2019faste} for more advantages of block methods. Numerical examples are given to illustrate the efficiency of REBK.

{\it Organization of the paper}. In the rest of this section, we give some notation. In Section 2 we describe the randomized extended block Kaczmarz algorithm and prove its convergence theory. Both the exponential convergence of the norm of the expected error and the exponential convergence of the expected norm of the error are discussed. In Section 3 we report the numerical results. Finally, we present brief concluding remarks in Section 4.

{\it Notation}. For any random variable $\bm\xi$, let $\mbbe\bem\bm\xi\eem$ denote its expectation. For an integer $m\geq 1$, let $[m]:=\{1,2,3,\ldots,m\}$. Lowercase (upper-case) boldface letters are reserved for column vectors (matrices). For any vector $\mbf u\in\mbbr^m$, we use $\mbf u_i$, $\bf u^\rmt$, and $\|\mbf u\|_2$ to denote the $i$th element, the transpose, and the $\ell_2$-norm of $\mbf u$, respectively. We use $\mbf I$ to denote the identity matrix whose order is clear from the context. For any matrix $\mbf A\in\mbbr^{m\times n}$, we use $\mbf A^\rmt$, $\mbf A^\dag$, $\|\mbf A\|_\rmf$, $\ran(\mbf A)$, $\sigma_{1}(\mbf A)\geq\sigma_{2}(\mbf A)\geq\cdots\geq\sigma_{r}(\mbf A)>0$ to denote the transpose, the Moore-Penrose pseudoinverse, the Frobenius norm, the column space, and all the nonzero singular values of $\mbf A$, respectively. Obviously, $r$ is the rank of $\mbf A$. For index sets $\mcali\subseteq[m]$ and $\mcalj\subseteq[n]$, let $\mbf A_{\mcali,:}$, $\mbf A_{:,\mcalj}$, and $\mbf A_{\mcali,\mcalj}$ denote the row submatrix indexed by $\mcali$, the column submatrix indexed by $\mcalj$, and the submatrix that lies in the rows indexed by $\mcali$ and the columns indexed by $\mcalj$, respectively. We call $\{\mcali_1,\mcali_2,\ldots,\mcali_s\}$ a partition of $[m]$ if $\mcali_i\cap\mcali_j=\emptyset$ for $i\neq j$ and $\cup_{i=1}^s\mcali_i=[m]$. Similarly, $\{\mcalj_1,\mcalj_2,\ldots,\mcalj_t\}$ denotes a partition of $[n]$ if $\mcalj_i\cap\mcalj_j=\emptyset$ for $i\neq j$ and $\cup_{j=1}^t\mcalj_j=[n]$. { We use $|\mcali|$ to denote the cardinality of a set $\mcali\subseteq[m]$}.

\section{The randomized extended block Kaczmarz algorithm} 
In this section, based on given partitions of $[m]$ and $[n]$, we propose the following  randomized extended block Kaczmarz algorithm (see Algorithm 1) for solving consistent or inconsistent linear systems. { We emphasize that this algorithm can be implemented on distributed computing units}. 

\begin{center}
\begin{tabular*}{150mm}{l}
\toprule {\bf Algorithm 1:} Randomized extended block Kaczmarz (REBK) \\ 
\hline
\qquad Let $\{\mcali_1,\mcali_2,\ldots,\mcali_s\}$ and $\{\mcalj_1,\mcalj_2,\ldots,\mcalj_t\}$ be partitions of $[m]$ and $[n]$, respectively.\\ 
\qquad Let $\alpha>0$. Initialize $\mbf z^0\in\mbbr^m$ and $\mbf x^0\in\mbbr^n$.\\
\qquad {\bf for} $k=1,2,\ldots,$ {\bf do}\\
\qquad\qquad Pick $j\in[t]$ with probability $\|{\bf A}_{:,\mcalj_j}\|_\rmf^2/\|{\bf A}\|_\rmf^2$\\
\qquad\qquad Set $\mbf z^k=\mbf z^{k-1}-\dsp\frac{\alpha}{\|{\bf A}_{:,\mcalj_j}\|_\rmf^2}\mbf A_{:,\mcalj_j}(\mbf A_{:,\mcalj_j})^\rmt\mbf z^{k-1}$\\
\qquad\qquad Pick $i\in[s]$ with probability $\|{\bf A}_{\mcali_i,:}\|_\rmf^2/\|{\bf A}\|_\rmf^2$\\
\qquad\qquad Set ${\bf x}^k={\bf x}^{k-1}-\dsp\frac{\alpha}{\|{\bf A}_{\mcali_i,:}\|_\rmf^2}(\mbf A_{\mcali_i,:})^\rmt(\mbf A_{\mcali_i,:}\mbf x^{k-1}-\mbf b_{\mcali_i}+\mbf z^k_{\mcali_i})$\\
\bottomrule
\end{tabular*}
\end{center}
Here we only consider constant stepsize for simplicity. By choosing the row partition parameter $s=m$, the column partition parameter $t=n$, and the stepsize $\alpha=1$, we recover the well-known randomized extended Kaczmarz algorithm of Zouzias and Freris \cite{zouzias2013rando}. REBK uses two RABK updates (see (\ref{rabk})) at each step: 
\bit
\item RABK update for $\bf A^\rmt z= 0$ from $\mbf z^{k-1}$
$$\mbf z^k=\mbf z^{k-1}-\alpha\l(\sum_{l\in\mcalj_j}\omega_l^k\frac{(\mbf A_{:,l})^\rmt\mbf z^{k-1}}{(\mbf A_{:,l})^\rmt\mbf A_{:,l}}\mbf A_{:,l}\r),\qquad \omega_l^k=\frac{(\mbf A_{:,l})^\rmt\mbf A_{:,l}}{\|\mbf A_{:,\mcalj_j}\|_\rmf^2};$$
\item RABK update for ${\bf Ax = b-z}^k$ from $\mbf x^{k-1}$
$$\mbf x^k=\mbf x^{k-1}-\alpha\l(\sum_{l\in\mcali_i}\omega_l^k\frac{\mbf A_{l,:}\mbf x^{k-1}-\mbf b_l+\mbf z_l^k}{\mbf A_{l,:}(\mbf A_{l,:})^\rmt}(\mbf A_{l,:})^\rmt\r), \qquad \omega_l^k=\frac{\mbf A_{l,:}(\mbf A_{l,:})^\rmt}{\|\mbf A_{\mcali_i,:}\|_\rmf^2}.$$
\eit  
{ We note that if $\mbf z^0=\mbf 0$ in REBK, then all $\mbf z^k\equiv\mbf 0$, which yields the update of $\mbf x^k$ is exactly the same as that of RABK.}

Before proving the convergence theory of REBK for general linear systems, we give the following notation. Let $\mbbe_{k-1}\bem\cdot\eem$ denote the conditional expectation conditioned on the first $k-1$ iterations of REBK. That is, $$\mbbe_{k-1}\bem\cdot\eem=\mbbe\bem\cdot|j_1,i_1,j_2,i_2,\ldots,j_{k-1},i_{k-1}\eem,$$ where $j_l$ is the $l$th column block chosen and $i_l$ is the $l$th row block chosen. We denote the conditional expectation conditioned on the first $k-1$ iterations and the $k$th column block chosen as $$\mbbe_{k-1}^i\bem\cdot\eem=\mbbe\bem\cdot|j_1,i_1,j_2,i_2,\ldots,j_{k-1},i_{k-1},j_k\eem.$$ 
Then by the law of total expectation we have $$\mbbe_{k-1}\bem\cdot\eem=\mbbe_{k-1}\bem\mbbe_{k-1}^i\bem\cdot\eem\eem.$$

\subsection{The exponential convergence of the norm of the expected error}
In this subsection we show the exponential convergence of the norm of the expected error, i.e., $${\|\mbbe\bem\mbf x^k-\mbf A^\dag\mbf b\eem\|_2}.$$ The convergence of the norm of the expected error depends on the positive number $\delta$ defined as $$\delta:=\max_{1\leq i\leq r}\l|1-\frac{\alpha\sigma_i^2(\mbf A)}{\|\mbf A\|_\rmf^2}\r|.$$ The following lemma will be used and its proof is straightforward (e.g., via the singular value decomposition). 

\begin{lemma}\label{leqd} Let $\alpha>0$ and $\mbf A\in\mbbr^{m\times n}$ be any nonzero real matrix with $\rank(\mbf A)=r$. For every $\mbf u\in\ran(\mbf A^\rmt)$, it holds $$\l\|\l(\mbf I-\alpha\frac{\bf A^\rmt A}{\|\mbf A\|_\rmf^2}\r)^k\mbf u\r\|_2\leq\delta^k\|\mbf u\|_2.$$ 
\end{lemma}

We give the convergence of the norm of the expected error of REBK in the following theorem.

\begin{theorem}\label{ner} For any given consistent or inconsistent linear system $\bf Ax=b$, let $\mbf x^k$ be the $k$th iterate of {\rm REBK} with {${\bf z}^0\in\mbbr^m$} and $\mbf x^0\in\ran(\mbf A^\rmt)$. It holds $$\|\mbbe\bem{\bf x}^k-{\bf A^\dag b}\eem\|_2\leq\delta^k\l(\|{\bf x}^0-{\bf A^\dag b}\|_2+\frac{\alpha k\|\mbf A^\rmt\mbf z^0\|_2}{\|\mbf A\|_\rmf^2}\r).$$  
\end{theorem} 
\proof Note that { \begin{align*} \mbbe_{k-1}\bem  \mbf z^k \eem &= \mbf z^{k-1} - \mbbe_{k-1}\bem \dsp\frac{\alpha}{\|{\bf A}_{:,\mcalj_j}\|_\rmf^2}\mbf A_{:,\mcalj_j}(\mbf A_{:,\mcalj_j})^\rmt \eem\mbf z^{k-1} \\ &= \l(\mbf I-\alpha\frac{\mbf A\mbf A^\rmt}{\|\mbf A\|_\rmf^2}\r)\mbf z^{k-1},\end{align*} and therefore $$\mbbe\bem \mbf z^k\eem=\mbbe\bem\mbbe_{k-1}\bem  \mbf z^k\eem \eem=\l(\mbf I-\alpha\frac{\mbf A\mbf A^\rmt}{\|\mbf A\|_\rmf^2}\r)\mbbe\bem \mbf z^{k-1}\eem=\l(\mbf I-\alpha\frac{\mbf A\mbf A^\rmt}{\|\mbf A\|_\rmf^2}\r)^k\mbf z^0.$$ }
By $\mbf A^\rmt\mbf b=\mbf A^\rmt\mbf A\mbf A^\dag\mbf b$, we have
\begin{align*}\mbbe_{k-1}\bem{\bf x}^k-{\bf A^\dag b}\eem
&=\mbbe_{k-1}\bem\mbbe_{k-1}^i\bem{\bf x}^k-{\bf A^\dag b}\eem\eem\\
&={\mbbe_{k-1}\bem\mbbe_{k-1}^i\bem{\bf x}^{k-1}-\mbf A^\dag \mbf b-\dsp\frac{\alpha}{\|{\bf A}_{\mcali_i,:}\|_\rmf^2}(\mbf A_{\mcali_i,:})^\rmt(\mbf A_{\mcali_i,:}\mbf x^{k-1}-\mbf b_{\mcali_i}+\mbf z^k_{\mcali_i})\eem\eem}\\
&={\mbbe_{k-1}\bem{\bf x}^{k-1}-\mbf A^\dag \mbf b- \alpha\dsp\frac{\mbf A^\rmt(\mbf A\mbf x^{k-1}-\mbf b+\mbf z^k)}{\|\mbf A\|_\rmf^2}\eem}\\
&= {\bf x}^{k-1}-\mbf A^\dag \mbf b-\alpha\dsp\frac{\mbf A^\rmt\mbf A\mbf x^{k-1}-\mbf A^\rmt\mbf b}{\|\mbf A\|_\rmf^2}- \alpha\dsp\frac{\mbf A^\rmt}{\|\mbf A\|_\rmf^2}\mbbe_{k-1}\bem\mbf z^k\eem\\
&={{\bf x}^{k-1}-\mbf A^\dag \mbf b-\alpha\dsp\frac{\mbf A^\rmt\mbf A\mbf x^{k-1}-\mbf A^\rmt\mbf A\mbf A^\dag\mbf b}{\|\mbf A\|_\rmf^2}- \alpha\dsp\frac{\mbf A^\rmt}{\|\mbf A\|_\rmf^2}\mbbe_{k-1}\bem\mbf z^k\eem}\\
&= \l(\mbf I-\alpha\frac{\mbf A^\rmt\mbf A}{\|\mbf A\|_\rmf^2}\r)(\mbf x^{k-1}-\mbf A^\dag\mbf b)- \alpha\dsp\frac{\mbf A^\rmt}{\|\mbf A\|_\rmf^2}\l(\mbf I-\alpha\frac{\mbf A\mbf A^\rmt}{\|\mbf A\|_\rmf^2}\r)\mbf z^{k-1}.\end{align*}
Taking expectation gives 
\begin{align*}
\mbbe\bem{\bf x}^k-{\bf A^\dag b}\eem&=\mbbe\bem\mbbe_{k-1}\bem{\bf x}^k-{\bf A^\dag b}\eem\eem\\
&=\l(\mbf I-\alpha\frac{\mbf A^\rmt\mbf A}{\|\mbf A\|_\rmf^2}\r)\mbbe\bem\mbf x^{k-1}-\mbf A^\dag\mbf b\eem-  \alpha\dsp\frac{\mbf A^\rmt}{\|\mbf A\|_\rmf^2}\l(\mbf I-\alpha\frac{\mbf A\mbf A^\rmt}{\|\mbf A\|_\rmf^2}\r)\mbbe\bem\mbf z^{k-1}\eem\\
&=\l(\mbf I-\alpha\frac{\mbf A^\rmt\mbf A}{\|\mbf A\|_\rmf^2}\r)\mbbe\bem\mbf x^{k-1}-\mbf A^\dag\mbf b\eem-  \alpha\dsp\frac{\mbf A^\rmt}{\|\mbf A\|_\rmf^2}\l(\mbf I-\alpha\frac{\mbf A\mbf A^\rmt}{\|\mbf A\|_\rmf^2}\r)^k\mbf z^0\\
&=\l(\mbf I-\alpha\frac{\mbf A^\rmt\mbf A}{\|\mbf A\|_\rmf^2}\r)\mbbe\bem\mbf x^{k-1}-\mbf A^\dag\mbf b\eem-  \alpha\dsp\l(\mbf I-\alpha\frac{\mbf A^\rmt\mbf A}{\|\mbf A\|_\rmf^2}\r)^k\frac{\mbf A^\rmt\mbf z^0}{\|\mbf A\|_\rmf^2}\\
&=\l(\mbf I-\alpha\frac{\mbf A^\rmt\mbf A}{\|\mbf A\|_\rmf^2}\r)^2\mbbe\bem\mbf x^{k-2}-\mbf A^\dag\mbf b\eem- 2\alpha\l(\mbf I-\alpha\frac{\mbf A^\rmt\mbf A}{\|\mbf A\|_\rmf^2}\r)^k\dsp\frac{\mbf A^\rmt\mbf z^0}{\|\mbf A\|_\rmf^2}\\&=\cdots\\
&=\l(\mbf I-\alpha\frac{\mbf A^\rmt\mbf A}{\|\mbf A\|_\rmf^2}\r)^k(\mbf x^0-\mbf A^\dag\mbf b)- \alpha k\l(\mbf I-\alpha\frac{\mbf A^\rmt\mbf A}{\|\mbf A\|_\rmf^2}\r)^k\dsp\frac{\mbf A^\rmt\mbf z^0}{\|\mbf A\|_\rmf^2}.\end{align*}
Applying the norms to both sides we obtain \begin{align*}\|\mbbe\bem{\bf x}^k-{\bf A^\dag b}\eem\|_2 &= \l\|\l(\mbf I-\alpha\frac{\mbf A^\rmt\mbf A}{\|\mbf A\|_\rmf^2}\r)^k(\mbf x^0-\mbf A^\dag\mbf b)- \alpha k\l(\mbf I-\alpha\frac{\mbf A^\rmt\mbf A}{\|\mbf A\|_\rmf^2}\r)^k\dsp\frac{\mbf A^\rmt\mbf z^0}{\|\mbf A\|_\rmf^2}\r\|_2\\ &\leq  \l\|\l(\mbf I-\alpha\frac{\mbf A^\rmt\mbf A}{\|\mbf A\|_\rmf^2}\r)^k(\mbf x^0-\mbf A^\dag\mbf b)\r\|_2+ \l\|\alpha k\l(\mbf I-\alpha\frac{\mbf A^\rmt\mbf A}{\|\mbf A\|_\rmf^2}\r)^k\dsp\frac{\mbf A^\rmt\mbf z^0}{\|\mbf A\|_\rmf^2}\r\|_2
\\&\leq \delta^k\l(\|{\bf x}^0-{\bf A^\dag b}\|_2+\frac{\alpha k\|\mbf A^\rmt\mbf z^0\|_2}{\|\mbf A\|_\rmf^2}\r).\end{align*}
Here the last inequality follows from the fact that {$\mbf x^0\in\ran(\mbf A^\rmt)$, ${\bf A^\dag b}\in\ran(\mbf A^\rmt)$}, $\mbf A^\rmt\mbf z^0\in\ran(\mbf A^\rmt)$, and Lemma \ref{leqd}.
\qed

\begin{remark}\label{alpha} To ensure convergence of the expected error, it suffices to have $$\delta=\dsp\max_{1\leq i\leq r}\l|1-\frac{\alpha\sigma_i^2(\mbf A)}{\|\mbf A\|_\rmf^2}\r|<1,$$ which implies $$0<\alpha<\dsp\frac{2\|\mbf A\|_\rmf^2}{\sigma_1^2(\mbf A)}.$$ {The optimal $\alpha$ in Theorem  \ref{ner} is  (see \cite{poljak1963gradi})  $$\frac{2\|\mbf A\|_\rmf^2}{\sigma_1^2(\mbf A)+\sigma_r^2(\mbf A)}=\argmin_{0<\alpha<\frac{2\|\mbf A\|_\rmf^2}{\sigma_1^2(\mbf A)}}\ \max_{1\leq i\leq r}\l|1-\frac{\alpha\sigma_i^2(\mbf A)}{\|\mbf A\|_\rmf^2}\r| ,$$ and the corresponding convergence rate $\delta$ is $$\frac{\sigma_1^2(\mbf A)-\sigma_r^2(\mbf A)}{\sigma_1^2(\mbf A)+\sigma_r^2(\mbf A)}.$$} 
\end{remark}

\subsection{The exponential convergence of the expected norm of the error}

In this subsection we show the exponential convergence of the expected norm of the error, i.e., $$\mbbe\bem\|\mbf x^k-\mbf A^\dag\mbf b\|_2^2\eem.$$ The convergence of the expected norm of the error depends on the positive numbers {$\eta$ and $\rho$ defined as $$\eta:=1-\frac{(2\alpha-\alpha^2\beta_{\max}^\mcali)\sigma_{r}^2(\bf A)}{\|\mbf A\|_\rmf^2},\qquad \rho:=1-\frac{(2\alpha-\alpha^2\beta_{\max}^\mcalj)\sigma_{r}^2(\bf A)}{\|\mbf A\|_\rmf^2},$$ where $$\beta_{\max}^\mcali:=\max_{i\in[s]}\frac{\|\mbf A_{\mcali_i,:}\|_2^2}{\|\mbf A_{\mcali_i,:}\|_\rmf^2},\qquad \beta_{\max}^\mcalj:=\max_{j\in[t]}\frac{\|\mbf A_{:,\mcalj_j}\|_2^2}{\|\mbf A_{:,\mcalj_j}\|_\rmf^2}.$$}The following lemmas will be used extensively in this paper.  
\begin{lemma}\label{leq} Let $\mbf A\in\mbbr^{m\times n}$ be any nonzero real matrix with $\rank(\mbf A)=r$. For every $\mbf u\in\ran(\mbf A)$, it holds $$\|{\bf A^\rmt}\mbf u\|_2^2\geq\sigma_r^2(\mbf A)\|\mbf u\|_2^2.$$
\end{lemma}
{\begin{lemma}\label{leqQ} Let $\mbf A\in\mbbr^{m\times n}$ be any nonzero real matrix. For every $\mbf u\in\mbbr^m$, it holds $$\mbf u^\rmt(\mbf A\mbf A^\rmt)^2 \mbf u\leq\|\mbf A\|_2^2\|\mbf A^\rmt\mbf u\|_2^2.$$ 
\end{lemma}
The proof of Lemma \ref{leq} is straightforward (e.g., via the singular value decomposition), and Lemma \ref{leqQ} follows from $$\mbf u^\rmt(\mbf A\mbf A^\rmt)^2 \mbf u=\mbf u^\rmt\mbf A(\mbf A^\rmt\mbf A)\mbf A^\rmt \mbf u\leq \|\mbf A^\rmt\mbf A\|_2\|\mbf A^\rmt\mbf u\|_2^2=\|\mbf A\|_2^2\|\mbf A^\rmt\mbf u\|_2^2.$$}In the following lemma we show that the vector $\mbf z^k$ generated in {\rm REBK} with  ${\bf z}^0\in{\bf b}+\ran({\bf A})$ converges to $$\mbf b_\bot=:(\bf I-AA^\dag)b,$$ which is the orthogonal projection of $\mbf z^0$ onto the set $\{\mbf z\ |\ \bf A^\rmt z=0\}$.
\begin{lemma}\label{estz} For any given consistent or inconsistent linear system $\bf Ax=b$, let $\mbf z^k$ be the vector generated in {\rm REBK} with ${\bf z}^0\in{\bf b}+\ran({\bf A})$. {Assume $0<\alpha<2/\beta_{\max}^\mcalj$.} It holds \beqs\mbbe\bem\|{\bf z}^k-\mbf b_\bot\|_2^2\eem\leq\rho^k\|\mbf z^0- \mbf b_\bot\|_2^2.\eeqs
\end{lemma}
\proof By $(\mbf A_{:,\mcalj_j})^\rmt \mbf b_\bot=\mbf 0$, we have \beq
{\bf z}^k-\mbf b_\bot={\bf z}^{k-1}-\mbf b_\bot-\dsp\frac{\alpha}{\|{\bf A}_{:,\mcalj_j}\|_\rmf^2}\mbf A_{:,\mcalj_j}(\mbf A_{:,\mcalj_j})^\rmt(\mbf z^{k-1}-\mbf b_\bot).\label{zkup}\eeq By ${\bf z}^0-\mbf b_\bot=\mbf A\mbf A^\dag\mbf z^0\in\ran({\bf A})$ { and $\mbf A_{:,\mcalj_j}(\mbf A_{:,\mcalj_j})^\rmt(\mbf z^{k-1}-\mbf b_\bot)\in\ran(\mbf A)$}, we can show that ${\bf z}^k-\mbf b_\bot\in\ran(\bf A)$ by induction.
It follows from (\ref{zkup}) that \begin{align*}\|{\bf z}^k-\mbf b_\bot\|_2^2
&=\|{\bf z}^{k-1}-\mbf b_\bot\|_2^2-\frac{2\alpha\|(\mbf A_{:,\mcalj_j})^\rmt(\mbf z^{k-1}-\mbf b_\bot)\|_2^2}{\|\mbf A_{:,\mcalj_j}\|_\rmf^2}\\
& \quad +\alpha^2(\mbf z^{k-1}-\mbf b_\bot)^\rmt\l(\frac{\mbf A_{:,\mcalj_j}}{\|{\bf A}_{:,\mcalj_j}\|_\rmf}\l(\frac{\mbf A_{:,\mcalj_j}}{\|{\bf A}_{:,\mcalj_j}\|_\rmf}\r)^\rmt\r)^2(\mbf z^{k-1}-\mbf b_\bot)\\
&\leq { \|{\bf z}^{k-1}-\mbf b_\bot\|_2^2-\l(2\alpha-\alpha^2\frac{\|\mbf A_{:,\mcalj_j}\|_2^2}{\|\mbf A_{:,\mcalj_j}\|_\rmf^2}\r)\frac{\|(\mbf A_{:,\mcalj_j})^\rmt(\mbf z^{k-1}-\mbf b_\bot)\|_2^2}{\|\mbf A_{:,\mcalj_j}\|_\rmf^2}}\\ & \qquad (\rm by\ Lemma\ \ref{leqQ})\\
&\leq  {\|{\bf z}^{k-1}-\mbf b_\bot\|_2^2-(2\alpha-\alpha^2\beta_{\max}^\mcalj)\frac{\|(\mbf A_{:,\mcalj_j})^\rmt(\mbf z^{k-1}-\mbf b_\bot)\|_2^2}{\|\mbf A_{:,\mcalj_j}\|_\rmf^2}}.\end{align*}
Taking the conditioned expectation on the first $k-1$ iterations   yields
\begin{align*}\mbbe_{k-1}\bem\|{\bf z}^k-\mbf b_\bot\|_2^2\eem &\leq  \|{\bf z}^{k-1}-\mbf b_\bot\|_2^2-\frac{(2\alpha-\alpha^2\beta_{\max}^\mcalj)\|\mbf A^\rmt(\mbf z^{k-1}-\mbf b_\bot)\|_2^2}{\|\mbf A \|_\rmf^2}\\
&\leq {\|{\bf z}^{k-1}-\mbf b_\bot\|_2^2-\frac{(2\alpha-\alpha^2\beta_{\max}^\mcalj)\sigma_r^2(\mbf A)}{\|\mbf A \|_\rmf^2}\|\mbf z^{k-1}-\mbf b_\bot\|_2^2} \\ &  {\qquad(\rm by\ Lemma\ \ref{leq}\ and\ 0<\alpha<2/\beta_{\max}^\mcalj)}\\
&= \rho\|{\bf z}^{k-1}-\mbf b_\bot\|_2^2 \end{align*} 
Taking expectation again gives
\begin{align*}\mbbe\bem\|{\bf z}^k-\mbf b_\bot\|_2^2\eem &=\mbbe\bem\mbbe_{k-1}\bem\|{\bf z}^k-\mbf b_\bot\|_2^2 \eem\eem\\
&\leq \rho\mbbe\bem\|{\bf z}^{k-1}-\mbf b_\bot\|_2^2\eem\\
&\leq \rho^k\|{\bf z}^0-\mbf b_\bot\|_2^2.\end{align*}
This completes the proof.\qed

We give the main convergence result of REBK in the following theorem.

\begin{theorem}\label{main} For any given consistent or inconsistent linear system $\bf Ax=b$, let $\mbf x^k$ be the $k$th iterate of {\rm REBK} with  ${\bf z}^0\in{\bf b}+\ran({\bf A})$ and $\mbf x^0\in\ran(\mbf A^\rmt)$. {Assume that $0<\alpha<2/\max(\beta_{\max}^\mcali,\beta_{\max}^\mcalj)$}. For any $\ve>0$, it holds $${\mbbe\bem\|{\bf x}^k-{\bf A^\dag b}\|_2^2\eem\leq(1+\ve)^k\eta^k\|{\bf x}^0-{\bf A^\dag b}\|_2^2+\l(1+\frac{1}{\ve}\r)\frac{\alpha^2\beta_{\max}^\mcali\|\mbf z^0-\mbf b_\bot\|_2^2}{\|\mbf A\|_\rmf^2}\sum_{l=0}^{k-1}\rho^{k-l}(1+\ve)^l\eta^l.}$$ 
\end{theorem}

\proof Let $$\wh{\mbf x}^k={\bf x}^{k-1}-\dsp\frac{\alpha}{\|{\bf A}_{\mcali_i,:}\|_\rmf^2}(\mbf A_{\mcali_i,:})^\rmt\mbf A_{\mcali_i,:}(\mbf x^{k-1}-\mbf A^\dag\mbf b),$$ which is actually one RABK update for the linear system $\bf Ax = AA^\dag b$ from $\mbf x^{k-1}$. It follows from  
$${\mbf x}^k-\wh{\mbf x}^k=\frac{\alpha}{\|\mbf A_{\mcali_i,:}\|_\rmf^2}(\mbf A_{\mcali_i,:})^\rmt(\mbf b_{\mcali_i}-\mbf A_{\mcali_i,:}\mbf A^\dag\mbf b-\mbf z^k_{\mcali_i})$$ that 
\begin{align}\|{\mbf x}^k-\wh{\mbf x}^k\|_2^2&=\frac{\alpha^2}{\|\mbf A_{\mcali_i,:}\|_\rmf^4}\|(\mbf A_{\mcali_i,:})^\rmt(\mbf b_{\mcali_i}-\mbf A_{\mcali_i,:}\mbf A^\dag\mbf b-\mbf z^k_{\mcali_i})\|_2^2\nn\\ 
&\leq  {\frac{\alpha^2}{\|\mbf A_{\mcali_i,:}\|_\rmf^2}\frac{\|\mbf A_{\mcali_i,:}\|_2^2}{\|\mbf A_{\mcali_i,:}\|_\rmf^2}\|\mbf b_{\mcali_i}-\mbf A_{\mcali_i,:}\mbf A^\dag\mbf b-\mbf z^k_{\mcali_i}\|_2^2}\nn\\ 
&\leq {\frac{\alpha^2\beta_{\max}^\mcali}{\|\mbf A_{\mcali_i,:}\|_\rmf^2}\|\mbf b_{\mcali_i}-\mbf A_{\mcali_i,:}\mbf A^\dag\mbf b-\mbf z^k_{\mcali_i}\|_2^2}. \label{xhat1} \end{align}
It follows from 
\begin{align*} \mbbe_{k-1}\bem \|{\mbf x}^k-\wh{\mbf x}^k\|_2^2\eem &= \mbbe_{k-1}\bem\mbbe_{k-1}^i\bem\|{\mbf x}^k-\wh{\mbf x}^k\|_2^2\eem\eem\\
& \leq {\mbbe_{k-1}\bem\mbbe_{k-1}^i\bem\dsp\frac{\alpha^2\beta_{\max}^\mcali}{\|\mbf A_{\mcali_i,:}\|_\rmf^2}\|\mbf b_{\mcali_i}-\mbf A_{\mcali_i,:}\mbf A^\dag\mbf b-\mbf z^k_{\mcali_i}\|_2^2\eem\eem \quad   \quad(\mbox{by (\ref{xhat1})})}\\
&= {\mbbe_{k-1}\bem\dsp\frac{\alpha^2\beta_{\max}^\mcali\|\mbf b -\mbf A \mbf A^\dag\mbf b-\mbf z^k \|_2^2}{\|\mbf A\|_\rmf^2}\eem}
\end{align*} that
\begin{align} \mbbe\bem \|{\mbf x}^k-\wh{\mbf x}^k\|_2^2\eem &= {\mbbe\bem\mbbe_{k-1}\bem \|{\mbf x}^k-\wh{\mbf x}^k\|_2^2 \eem\eem}\nn\\ &\leq{\frac{\alpha^2\beta_{\max}^\mcali}{\|\mbf A\|_\rmf^2}\mbbe\bem\|\mbf b -\mbf A \mbf A^\dag\mbf b-\mbf z^k \|_2^2\eem}\nn\\ 
&\leq {\frac{\alpha^2\beta_{\max}^\mcali\rho^k}{\|\mbf A\|_\rmf^2}\|\mbf z^0-\mbf b_\bot \|_2^2. \quad ({\rm by\ Lemma\ \ref{estz}})}
\label{xhat}\end{align}
By $\mbf x^0 \in\ran(\mbf A^\rmt)$, { ${\bf A^\dag b}\in\ran(\mbf A^\rmt)$, $(\mbf A_{\mcali_i,:})^\rmt(\mbf A_{\mcali_i,:}\mbf x^{k-1}-\mbf b_{\mcali_i}+\mbf z^k_{\mcali_i})\in\ran(\mbf A^\rmt)$, and $${\bf x}^k-\mbf A^\dag\mbf b={\bf x}^{k-1}-\mbf A^\dag\mbf b-\dsp\frac{\alpha}{\|{\bf A}_{\mcali_i,:}\|_\rmf^2}(\mbf A_{\mcali_i,:})^\rmt(\mbf A_{\mcali_i,:}\mbf x^{k-1}-\mbf b_{\mcali_i}+\mbf z^k_{\mcali_i}),$$} we can show that $\mbf x^k-{\bf A^\dag b}\in\ran(\mbf A^\rmt)$ by induction. {By \begin{align*}\|\wh{\mbf x}^k-{\bf A^\dag b}\|_2^2&= \|\mbf x^{k-1}-\mbf A^\dag\mbf b\|_2^2-\frac{2\alpha\|\mbf A_{\mcali_i,:}(\mbf x^{k-1}-\mbf A^\dag\mbf b)\|_2^2}{\|\mbf A_{\mcali_i,:}\|_\rmf^2}\\
&\quad+\alpha^2(\mbf x^{k-1}-\mbf A^\dag\mbf b)^\rmt\l(\l(\frac{\mbf A_{\mcali_i,:}}{\|\mbf A_{\mcali_i,:}\|_\rmf}\r)^\rmt\frac{\mbf A_{\mcali_i,:}}{\|\mbf A_{\mcali_i,:}\|_\rmf}\r)^2(\mbf x^{k-1}-\mbf A^\dag\mbf b)\\
&\leq\|\mbf x^{k-1}-\mbf A^\dag\mbf b\|_2^2-\frac{(2\alpha-\alpha^2\beta_{\max}^\mcali)\|\mbf A_{\mcali_i,:}(\mbf x^{k-1}-\mbf A^\dag\mbf b)\|_2^2}{\|\mbf A_{\mcali_i,:}\|_\rmf^2},\quad (\mbox{by Lemma \ref{leqQ}})\end{align*}
we have \begin{align*}\mbbe_{k-1}\bem\|\wh{\mbf x}^k-{\bf A^\dag b}\|_2^2\eem&\leq \|\mbf x^{k-1}-\mbf A^\dag\mbf b\|_2^2-\frac{(2\alpha-\alpha^2\beta_{\max}^\mcali)\|\mbf A(\mbf x^{k-1}-\mbf A^\dag\mbf b)\|_2^2}{\|\mbf A\|_\rmf^2}\\&\leq \|\mbf x^{k-1}-\mbf A^\dag\mbf b\|_2^2-\frac{(2\alpha-\alpha^2\beta_{\max}^\mcali)\sigma_r^2(\mbf A)\|(\mbf x^{k-1}-\mbf A^\dag\mbf b)\|_2^2}{\|\mbf A\|_\rmf^2}\\ &\quad \quad (\mbox{by Lemma \ref{leq}} \mbox{ and } 0<\alpha<2/\beta_{\max}^\mcali) \\ &= \eta \|{\bf x}^{k-1}-{\bf A^\dag b}\|_2^2,\end{align*}
which yields \beq\label{la}\mbbe\bem\|\wh{\mbf x}^k-{\bf A^\dag b}\|_2^2\eem\leq\eta\mbbe\bem\|{\bf x}^{k-1}-{\bf A^\dag b}\|_2^2\eem.\eeq}
Note that for any $\ve>0$, we have \begin{align}\|\mbf x^k-\mbf A^\dag\mbf b\|_2^2&= \|\mbf x^k-\wh{\mbf x}^k+\wh{\mbf x}^k-\mbf A^\dag\mbf b\|_2^2 \nn \\ &\leq (\|\mbf x^k-\wh{\mbf x}^k\|_2+\|\wh{\mbf x}^k-\mbf A^\dag\mbf b\|_2)^2\nn \\&\leq  \|\mbf x^k-\wh{\mbf x}^k\|_2^2+\|\wh{\mbf x}^k-\mbf A^\dag\mbf b\|_2^2+2\|\mbf x^k-\wh{\mbf x}^k\|_2\|\wh{\mbf x}^k-\mbf A^\dag\mbf b\|_2\nn \\
&\leq\l(1+\frac{1}{\ve}\r)\|\mbf x^k-\wh{\mbf x}^k\|_2^2+(1+\ve)\|\wh{\mbf x}^k-\mbf A^\dag\mbf b\|_2^2.\label{ve}\end{align} 
Combining (\ref{xhat}), (\ref{la}), and (\ref{ve}) yields {\begin{align*}\mbbe\bem\|{\mbf x}^k-{\bf A^\dag b}\|_2^2\eem&\leq \l(1+\frac{1}{\ve}\r)\mbbe\bem\|{\mbf x}^k-\wh{\mbf x}^k\|_2^2\eem+(1+\ve)\mbbe\bem\|\wh{\mbf x}^k-{\bf A^\dag b}\|_2^2\eem\\
&\leq \l(1+\frac{1}{\ve}\r)\frac{\alpha^2\beta_{\max}^\mcali\rho^k}{\|\mbf A\|_\rmf^2}\|\mbf z^0-\mbf b_\bot\|_2^2 +(1+\ve)\eta\mbbe\bem\|{\bf x}^{k-1}-{\bf A^\dag b}\|_2^2\eem\\
&\leq \l(1+\frac{1}{\ve}\r)\frac{\alpha^2\beta_{\max}^\mcali\|\mbf z^0-\mbf b_\bot\|_2^2}{\|\mbf A\|_\rmf^2}(\rho^k+\rho^{k-1}(1+\ve)\eta)\\ &\quad  +(1+\ve)^2\eta^2\mbbe\bem\|{\bf x}^{k-2}-{\bf A^\dag b}\|_2^2\eem\\
&\leq \cdots\\ &\leq \l(1+\frac{1}{\ve}\r)\frac{\alpha^2\beta_{\max}^\mcali\|\mbf z^0-\mbf b_\bot\|_2^2}{\|\mbf A\|_\rmf^2}\sum_{l=0}^{k-1}\rho^{k-l}(1+\ve)^l\eta^l\\ &\quad +(1+\ve)^k\eta^k\|{\bf x}^0-{\bf A^\dag b}\|_2^2.
\end{align*}}This completes the proof.
\qed

\begin{remark}\label{rek8} For the case REBK with $s=m$, $t=n$ and $\alpha=1$ (i.e., REK), we have {$$\beta_{\max}^\mcali=\max_{i\in[m]}\frac{\|\mbf A_{i,:}\|_2^2}{\|\mbf A_{i,:}\|_\rmf^2}=1,\quad \beta_{\max}^\mcalj=\max_{j\in[n]}\frac{\|\mbf A_{:,j}\|_2^2}{\|\mbf A_{:,j}\|_\rmf^2}=1.$$ Therefore,
$$\eta=1-\frac{(2\alpha-\alpha^2\beta_{\max}^\mcali)\sigma_{r}^2(\bf A)}{\|\mbf A\|_\rmf^2}=1-\frac{\sigma_r^2(\mbf A)}{\|\mbf A\|_\rmf^2},$$ and $$\rho=1-\frac{(2\alpha-\alpha^2\beta_{\max}^\mcalj)\sigma_{r}^2(\bf A)}{\|\mbf A\|_\rmf^2}=1-\frac{\sigma_r^2(\mbf A)}{\|\mbf A\|_\rmf^2}.$$ It follows from $$\wh{\mbf x}^k-\mbf A^\dag\mbf b=\l(\mbf I-\frac{(\mbf A_{i,:})^\rmt\mbf A_{i,:}}{\|\mbf A_{i,:}\|_2^2}\r)(\mbf x^{k-1}-\mbf A^\dag\mbf b)$$ and $${\mbf x}^k-\wh{\mbf x}^k=\frac{\mbf b_i-\mbf A_{i,:}\mbf A^\dag\mbf b-\mbf z^k_i}{\|\mbf A_{i,:}\|_2^2}(\mbf A_{i,:})^\rmt$$ that $$(\wh{\mbf x}^k-\mbf A^\dag\mbf b)^\rmt(\mbf x^k-\wh{\mbf x}^k)= 0.$$ Then we have $$\|\mbf x^k-\mbf A^\dag\mbf b\|_2^2=\|\mbf x^k-\wh{\mbf x}^k\|_2^2+\|\wh{\mbf x}^k-\mbf A^\dag\mbf b\|_2^2,$$ which yields the following convergence for REK (see \cite{du2019tight}):  \begin{align*}\mbbe\bem\|{\mbf x}^k-{\bf A^\dag b}\|_2^2\eem&= \mbbe\bem\|{\mbf x}^k-\wh{\mbf x}^k\|_2^2\eem+\mbbe\bem\|\wh{\mbf x}^k-{\bf A^\dag b}\|_2^2\eem\\
&\leq \frac{\alpha^2\rho^k}{\|\mbf A\|_\rmf^2}\|\mbf z^0-\mbf b_\bot\|_2^2 +\rho\mbbe\bem\|{\bf x}^{k-1}-{\bf A^\dag b}\|_2^2\eem\\
&\leq \frac{2\alpha^2\rho^k\|\mbf z^0-\mbf b_\bot\|_2^2}{\|\mbf A\|_\rmf^2}+\rho^2\mbbe\bem\|{\bf x}^{k-2}-{\bf A^\dag b}\|_2^2\eem\\
&\leq \cdots\\ &\leq \rho^k\l(\frac{k\|\mbf z^0-\mbf b_\bot\|_2^2}{\|\mbf A\|_\rmf^2}+\|{\bf x}^0-{\bf A^\dag b}\|_2^2\r).
\end{align*}
Actually our proof is a modification of that of Zouzias and Freris \cite{zouzias2013rando}. We reorganize the arguments used by Zouzias and Freris and refine the analysis to get a better convergence estimate.}
\end{remark}

\begin{remark}\label{rek9} {Let $
\wh\rho:=\max(\eta,\rho)$ and $\beta_{\max}:=\max(\beta_{\max}^\mcali,\beta_{\max}^\mcalj)$. Then we have $$\wh\rho=1-\frac{(2\alpha-\alpha^2\beta_{\max})\sigma_{r}^2(\bf A)}{\|\mbf A\|_\rmf^2}.$$ By Theorem \ref{main}, we have 
\begin{align*}\mbbe\bem\|{\mbf x}^k-{\bf A^\dag b}\|_2^2\eem
&\leq(1+\ve)^k\eta^k\|{\bf x}^0-{\bf A^\dag b}\|_2^2+\l(1+\frac{1}{\ve}\r)\frac{\alpha^2\beta_{\max}^\mcali\|\mbf z^0-\mbf b_\bot\|_2^2}{\|\mbf A\|_\rmf^2}\sum_{l=0}^{k-1}\rho^{k-l}(1+\ve)^l\eta^l\\
&\leq(1+\ve)^k\wh\rho^k\|{\bf x}^0-{\bf A^\dag b}\|_2^2+\l(1+\frac{1}{\ve}\r)\frac{\alpha^2\beta_{\max}^\mcali\|\mbf z^0-\mbf b_\bot\|_2^2}{\|\mbf A\|_\rmf^2}\wh\rho^k\sum_{l=0}^{k-1}(1+\ve)^l\\
&\leq(1+\ve)^k\wh\rho^k\|{\bf x}^0-{\bf A^\dag b}\|_2^2+\l(1+\frac{1}{\ve}\r)\frac{\alpha^2\beta_{\max}^\mcali\|\mbf z^0-\mbf b_\bot\|_2^2}{\|\mbf A\|_\rmf^2}\wh\rho^k\frac{(1+\ve)^k-1}{\ve}\\
&\leq (1+\ve)^k\wh\rho^k\l(\|{\bf x}^0-{\bf A^\dag b}\|_2^2+\frac{(1+\ve)\alpha^2\beta_{\max}^\mcali\|\mbf z^0-\mbf b_\bot\|_2^2}{\ve^2\|\mbf A\|_\rmf^2}\r),
\end{align*}
which shows that REBK exponentially converges in the mean square to the minimum $\ell_2$-norm least squares solution of a given linear system of equations with the rate $(1+\ve)\wh\rho$  if  $\dsp 0<\alpha<{2}/{\beta_{\max}}$. Setting $\alpha=1/\beta_{\max}$ yields $$\wh\rho=1-\frac{\sigma_{r}^2(\bf A)}{\beta_{\max}\|\mbf A\|_\rmf^2},$$ which is better than the rate of REK (see Remark \ref{rek8}) $$\rho=\dsp1-\frac{\sigma_{r}^2(\bf A)}{\|\mbf A\|_\rmf^2}$$ if $\beta_{\max}<1$. A smaller $\beta_{\max}$ means a faster convergence in terms of iterations. Recalling that $$\beta_{\max}^\mcali:=\max_{i\in[s]}\frac{\|\mbf A_{\mcali_i,:}\|_2^2}{\|\mbf A_{\mcali_i,:}\|_\rmf^2}\quad \mbox{ and }\quad \beta_{\max}^\mcalj:=\max_{j\in[t]}\frac{\|\mbf A_{:,\mcalj_j}\|_2^2}{\|\mbf A_{:,\mcalj_j}\|_\rmf^2},$$ we have $$\max_{i\in[s]}\frac{1}{|\mcali_i|}\leq\max_{i\in[s]}\frac{1}{\rank(\mbf A_{\mcali_i,:})}\leq\beta_{\max}^\mcali\leq 1$$ and $$\max_{j\in[t]}\frac{1}{|\mcalj_j|}\leq\max_{j\in[t]}\frac{1}{\rank(\mbf A_{:,\mcali_j})}\leq\beta_{\max}^\mcalj\leq 1.$$ Therefore, $$\max\l(\max_{i\in[s]}\frac{1}{|\mcali_i|},\max_{j\in[t]}\frac{1}{|\mcalj_j|}\r)\leq\beta_{\max} \leq 1,$$ which means that REBK is at least as fast as REK in terms of iterations. The numerical results in Section 3 show that the convergence of REBK with appropriate block sizes and stepsizes is much faster than that of REK both in the numbers of iterations and the computing times.} 
\end{remark}

\begin{remark} It was shown in \cite{gower2015rando} that the convergence of $\mbf x^k$ to $\mbf A^\dag\mbf b$ under the expected norm of the error (Theorem \ref{main}) is a stronger form of convergence than the convergence of the norm of the expected error (Theorem \ref{ner}), as the former also guarantees that the variance of $\mbf x^k_i$ (the $i$th element of $\mbf x^k$) converges to zero for $i=1,\ldots,n$. By Remark \ref{alpha}, we know $0<\alpha<{2\|\mbf A\|_\rmf^2}/{\sigma_1^2(\mbf A)}$ guarantees the convergence of the norm of the expected error. {By Remark \ref{rek9}, we know $0<\alpha<2/\beta_{\max}$ guarantees the convergence of the expected norm of the error. However, since the convergence estimate in Remark \ref{rek9} usually is not sharp, the stepsize $\alpha$ satisfying $2/\beta_{\max}\leq\alpha<{2\|\mbf A\|_\rmf^2}/{\sigma_1^2(\mbf A)}$ is also possible to result in convergence (see Figure \ref{fig1}, Tables \ref{t2} and \ref{t3} in Section 3).}   
\end{remark}

\section{Numerical results}
In this section, we compare the performance of the randomized extended block Kaczmarz (REBK) algorithm proposed in this paper against the randomized extended Kaczmarz (REK) algorithm \cite{zouzias2013rando} and the projection-based randomized double block Kaczmarz (RDBK) algorithm \cite{needell2015rando}  on a variety of test problems. We do not claim optimized implementations of the algorithms, and only run on small or medium-scale problems. The purpose is only to demonstrate that even in these simple examples, REBK offers significant advantages to REK. All experiments are performed using MATLAB (version R2019a) on a laptop with 2.7-GHz Intel Core i7 processor, 16 GB memory, and Mac operating system. 

To construct an inconsistent linear system, we set $\bf b = Ax+r$ where $\mbf x$ is a vector with entries generated from a standard normal distribution and the residual $\mbf r\in \nul(\mbf A^\rmt)$. Note that one can obtain such a vector $\bf r$ by the MATLAB function $\tt null$. For all algorithms, we set $\mbf z^0=\mbf b$ and $\mbf x^0=\mbf 0$ and stop if the error $\|\mbf x^k-\mbf A^\dag\mbf b\|_2\leq10^{-5}$. 
We report the average number of iterations (denoted as ITER) {and the average computing time in seconds (denoted as CPU)}  of REK, RDBK, and REBK. { Note that {\tt A$\backslash$b} will usually not be the same as {\tt pinv(A)*b} when {\tt A} is rank-deficient or underdetermined. We use MATLAB's {\tt lsqminnorm} (which is typically more efficient than {\tt pinv}) to solve the small least squares problems at each step of RDBK.} We refer the reader to \cite{needell2014paved,needell2015rando} for more numerical aspects of RDBK. We also report the speed-up of REBK against REK, which is defined as $$\mbox{speed-up}=\frac{\mbox{CPU of REK}}{\mbox{CPU of REBK}}.$$  

For the block methods, we assume that the subsets $\{\mcali_i\}_{i=1}^{s-1}$ and $\{\mcalj_j\}_{j=1}^{t-1}$ have the same size $\tau$ (i.e., $|\mcali_i|=|\mcalj_j|=\tau$). We consider the row partition $\{\mcali_i\}_{i=1}^s$: \begin{align*}\mcali_i &=\{(i-1)\tau+1,(i-1)\tau+2,\ldots,i\tau\},\quad i=1,2,\ldots,s-1,\\ \mcali_s &=\{(s-1)\tau+1,(s-1)\tau+2,\ldots,m\},\quad |\mcali_s|\leq\tau,\end{align*} and the column partition $\{\mcalj_j\}_{j=1}^t$: \begin{align*}\mcalj_j &=\{(j-1)\tau+1,(j-1)\tau+2,\ldots,j\tau\},\quad j=1,2, \ldots, t-1,\\ \mcalj_t &= \{(t-1) \tau+1,(t-1)\tau+2, \ldots, n\},\quad |\mcalj_t|\leq \tau.\end{align*} 

\subsection{Synthetic data}

Two types of coefficient matrices are generated as follows. 
\bit
\item Type I: For given $m$, $n$, $r = rank(\mbf A)$, and $\kappa>1$, we construct a matrix $\mbf A$ by $$\bf A = UDV^\rmt,$$ where $\mbf U\in\mbbr^{m\times r} $ and $\mbf V\in \mbbr^{n\times r}$. Entries of $\mbf U$ and $\mbf V$ are generated from a standard normal distribution, and then, columns are orthonormalized, $${\tt [U,\sim] = qr(randn(m,r),0);\qquad  [V,\sim] = qr(randn(n,r),0);}$$ The matrix $\mbf D$ is an $r\times r$ diagonal matrix whose diagonal entries are uniformly distributed numbers in $(1,\kappa)$, $${\tt D = diag(1+(\kappa-1).*rand(r,1));}$$ So the condition number of $\mbf A$ is upper bounded by $\kappa$.
\item Type II: For given $m$, $n$, entries of $\mbf A$  are generated from a standard normal distribution, $$\tt A=randn(m,n);$$ So $\mbf A$ is a full-rank matrix almost surely.     
\eit

\begin{figure}[!htpb]
\centerline{\epsfig{figure=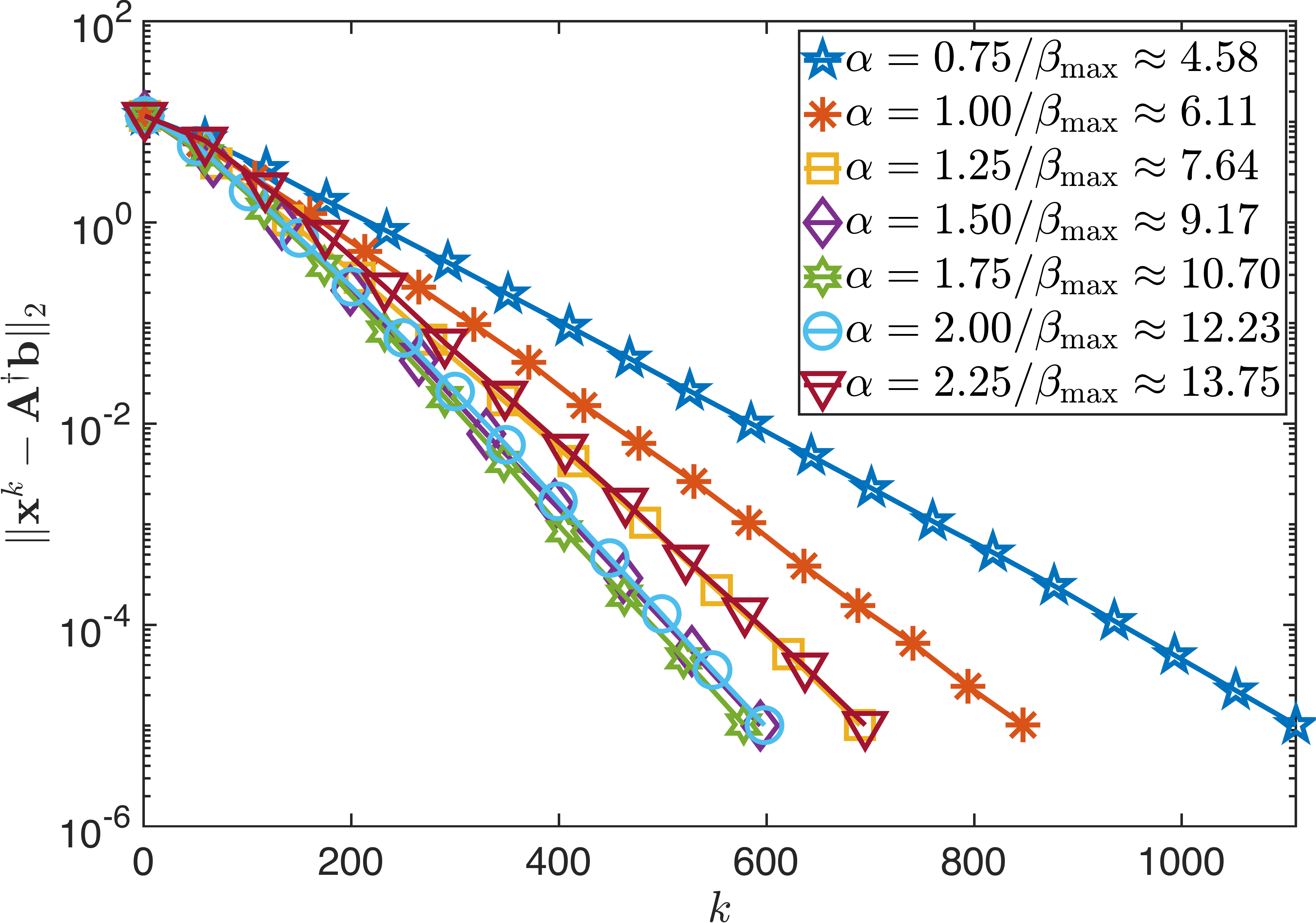,height=2.15in}\quad\epsfig{figure=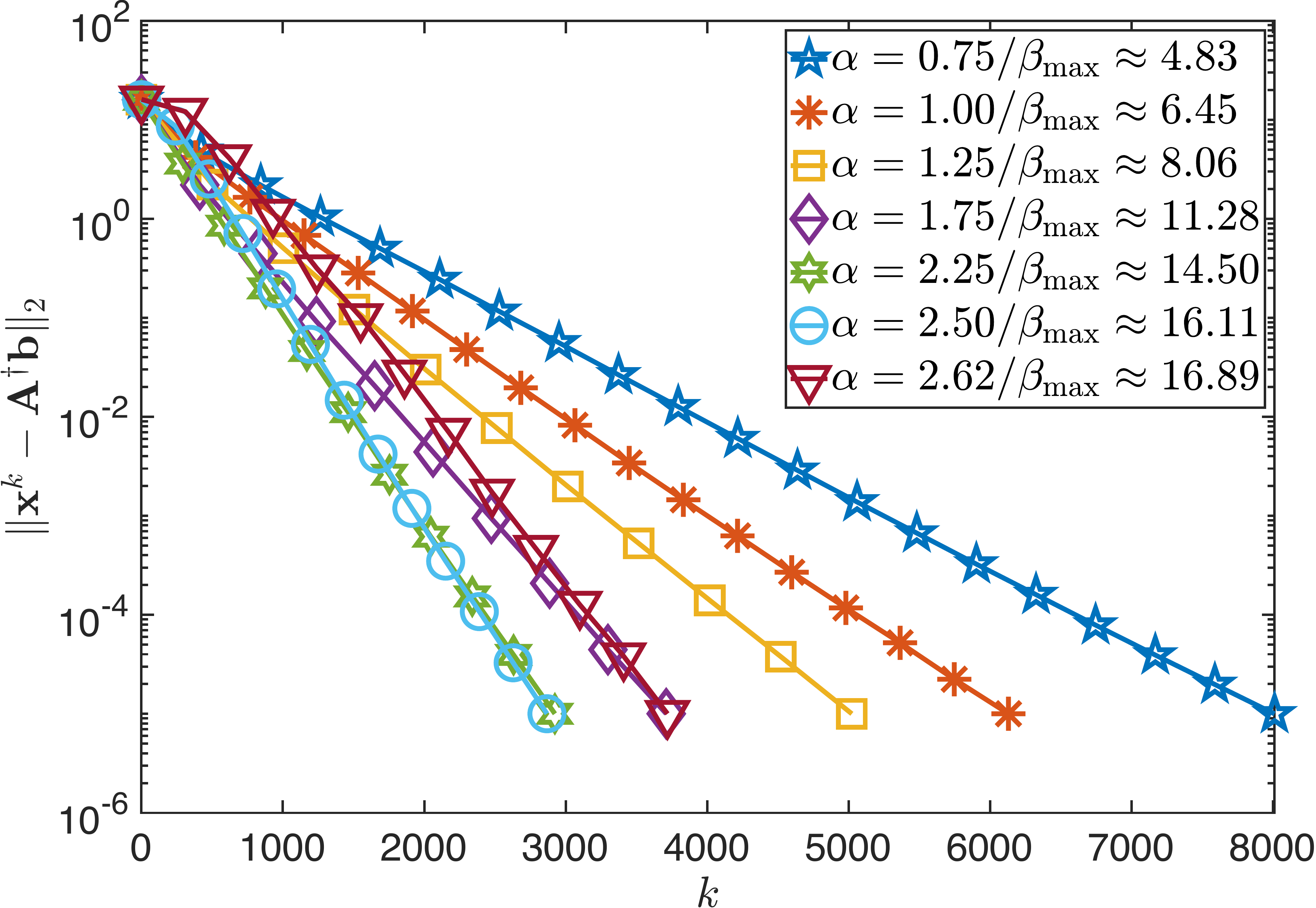,height=2.15in}}
\caption{The average (10 trials of each case) error $\|\mbf x^k-\mbf A^\dag\mbf b\|_2$ of REBK with block size $\tau=10$ and different stepsizes $\alpha$ from $0.75/\beta_{\max}$ to $2.62/\beta_{\max}$ for two inconsistent linear systems. Left: Type I matrix $\bf A=UDV^\rmt$ with $m=500$, $n=250$, $r=150$, $\kappa=2$. Right: Type II matrix $\bf A=$ {\tt randn(500,250)}.} \label{fig1} 
\end{figure}

{ In Figure \ref{fig1}, we plot the error $\|\mbf x^k-\mbf A^\dag \mbf b\|_2$ of REBK with a fixed block size ($\tau=10$) and different stepsizes ($\alpha$ from $0.75/\beta_{\max}$ to $2.62/\beta_{\max}$) for two inconsistent linear systems with coefficient matrices of Types I ($\bf A=UDV^\rmt$ with $m=500$, $n=250$, $r=150$, $\kappa=2$) and II ($\bf A=$ {\tt randn(500,250)}). It is observed that the convergence   of REBK becomes faster as the increase of the stepsize, and then slows down after reaching the fastest rate.}

{ In Tables \ref{t1} and \ref{t2}, we report the numbers of iterations and the computing times of the REK, RDBK, and REBK algorithms for solving inconsistent linear systems. For the block algorithms (RDBK and REBK), a fixed block size $\tau=10$ is used. For the REBK algorithm, empirical stepsizes $\alpha=1.75/\beta_{\max}$ and $\alpha=2.25/\beta_{\max}$ are used for Type I and Type II matrices, respectively. From these two tables, we observe: (i) in all cases, the RDBK and REBK algorithms vastly outperform the REK algorithm in terms of both the numbers of iterations and the computing times; (ii) for Type I matrix, the convergence rates of the RDBK and REBK algorithms are almost the same in terms of the numbers of iterations; (iii) for Type II matrix, REBK performs better than RDBK in terms of the numbers of iterations.}

\begin{table}[!htp]
\caption{The average (10 trials of each algorithm) {\rm ITER} and {\rm CPU} of {\rm REK}, {\rm RDBK}($\tau=10$), and {\rm REBK}($\tau=10$, $\alpha=1.75/\beta_{\max}$) for inconsistent linear systems  with  random coefficient matrices $\mbf A$ of Type I: ${\bf A=UDV^\rmt}$.} 
\label{t1}
\begin{center} \footnotesize
\begin{tabular}{c|c|c|c|c|c|c|c|c|c|c} \toprule
 \multirow{2}{*}{$m\times n$}& \multirow{2}{*}{rank} & \multirow{2}{*}{$\kappa$} & \multicolumn{2}{|c|}{REK}  & \multicolumn{2}{|c}{RDBK} & \multicolumn{4}{|c}{REBK}\\
  \cline{4-11}&&&ITER &CPU & ITER &CPU &
  $\alpha$ & ITER & CPU & speed-up\\
\hline
$ 250 \times 500 $ & 150 & 2 &5826 &0.26 &572 &0.21 &10.87 &586 &0.05 &4.90 \\ \hline
$ 250 \times 500 $ & 150 & 10 &65520 &2.87 &6166 &2.19 &9.36 &7365 &0.63 &4.59 \\ \hline
$ 500\times 1000$ & 250 & 2 &10068 &0.59 &1000 &0.43 &11.82 &991 &0.13 &4.60 \\ \hline
$ 500\times 1000$ & 250 & 10 &114297 &6.61 &10209 &4.29 &10.85 &10259 &1.23 &5.36 \\ \hline
$ 500\times 250$ & 150 & 2 &5755 &0.25 &562 &0.19 &10.70 &578 &0.03 &7.32 \\ \hline
$ 500\times 250$ & 150 & 10 &63741 &2.76 &5784 &1.90 &10.13 &6424 &0.36 &7.60 \\ \hline
$ 500\times 250$ & 250 & 2 &9971 &0.43 &940 &0.31 &12.47 &961 &0.06 &7.81 \\ \hline
$ 500\times 250$ & 250 & 10 &119182 &5.14 &11328 &3.73 &10.99 &10783 &0.61 &8.43 \\ \hline
$ 1000\times 500$ & 250 & 2 &9959 &0.55 &974 &0.39 &12.10 &987 &0.10 &5.53 \\ \hline
$ 1000\times 500$ & 250 & 10 &118134 &6.54 &11236 &4.44 &11.20 &10349 &1.03 &6.36\\ \hline
$ 1000\times 500$ & 500 & 2 &20188 &1.11 &2007 &0.80 &13.84 &2115 &0.21 &5.20 \\ \hline
$ 1000\times 500$ & 500 & 10 &254117 &14.01 &25361 &10.00 &12.67 &20432 &2.03 &6.92 \\
\bottomrule
\end{tabular}
\end{center}
\end{table}

\begin{table}[!htp]
\caption{The average (10 trials of each algorithm) {\rm ITER} and {\rm CPU} of {\rm REK}, {\rm RDBK}($\tau=10$), and {\rm REBK}($\tau=10$, $\alpha=2.25/\beta_{\max}$) for inconsistent linear systems  with  random coefficient matrices $\mbf A$ of Type II: ${\bf A=}$ {\tt randn(m,n)}.} 
\label{t2}
\begin{center} \footnotesize
\begin{tabular}{c|c|c|c|c|c|c|c|c|c|c} \toprule
 \multirow{2}{*}{$m\times n$}& \multirow{2}{*}{rank} & \multirow{2}{*}{$\dsp\frac{\sigma_1(\mbf A)}{\sigma_r(\mbf A)}$} & \multicolumn{2}{|c|}{REK} & \multicolumn{2}{|c}{RDBK} & \multicolumn{4}{|c}{REBK}\\
  \cline{4-11}&&&ITER &CPU &ITER &CPU &
  $\alpha$ & ITER & CPU & speed-up\\
\hline
$ 250\times 120$ & 120 &5.25  &18060 &0.66 &1646 &0.50 &13.48 &1337 &0.06 &10.54 \\ \hline
$ 500\times 250$ & 250 & 5.73 & 41016  &1.79  &3811 &1.26 &14.50 &2885 &0.17 &10.81 \\ \hline
$ 750\times 370$ & 370 & 5.80 &59660 &2.91 &5929 &2.20 &16.23 &4115 &0.36 &8.07 \\ \hline
$ 1000\times 500$ & 500 &5.74  &83093 &4.61 & 8183 & 3.25 & 16.42 &5422 &0.55 &8.41 \\ 
\bottomrule
\end{tabular}
\end{center}
\end{table}

\begin{figure}[htp]
\centerline{\epsfig{figure=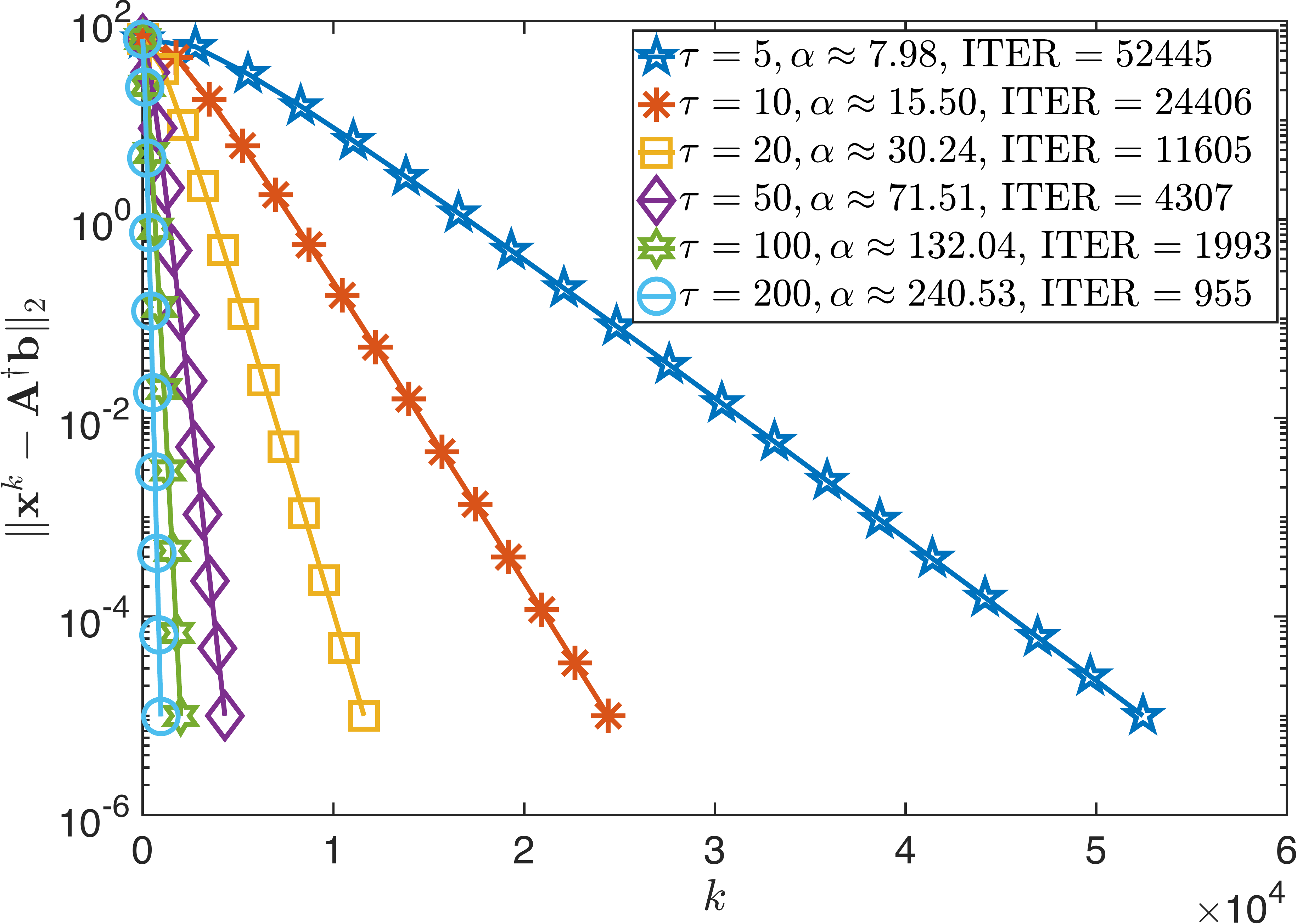,height=2.15in}\quad\epsfig{figure=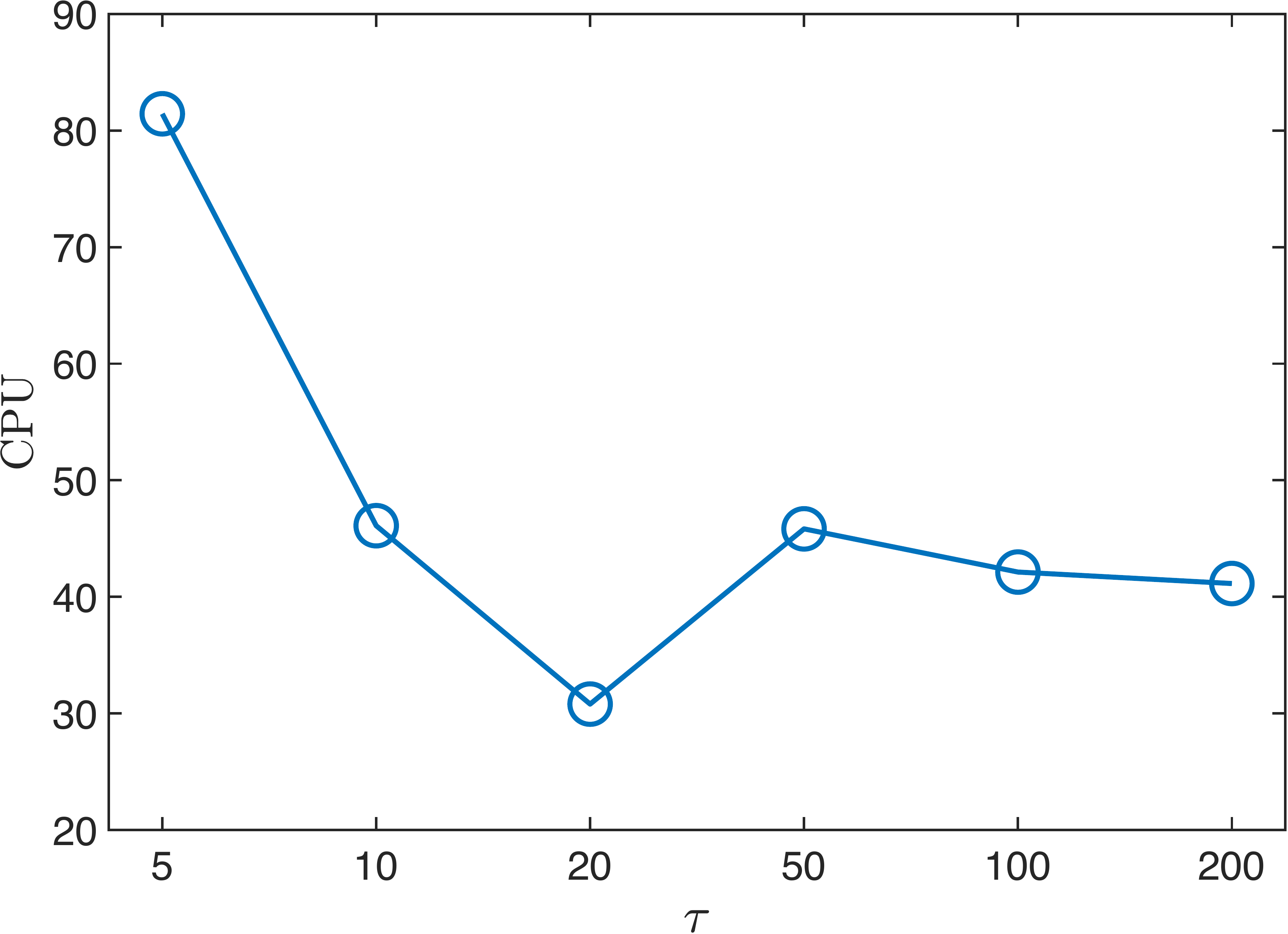,height=2.15in}}
\centerline{\epsfig{figure=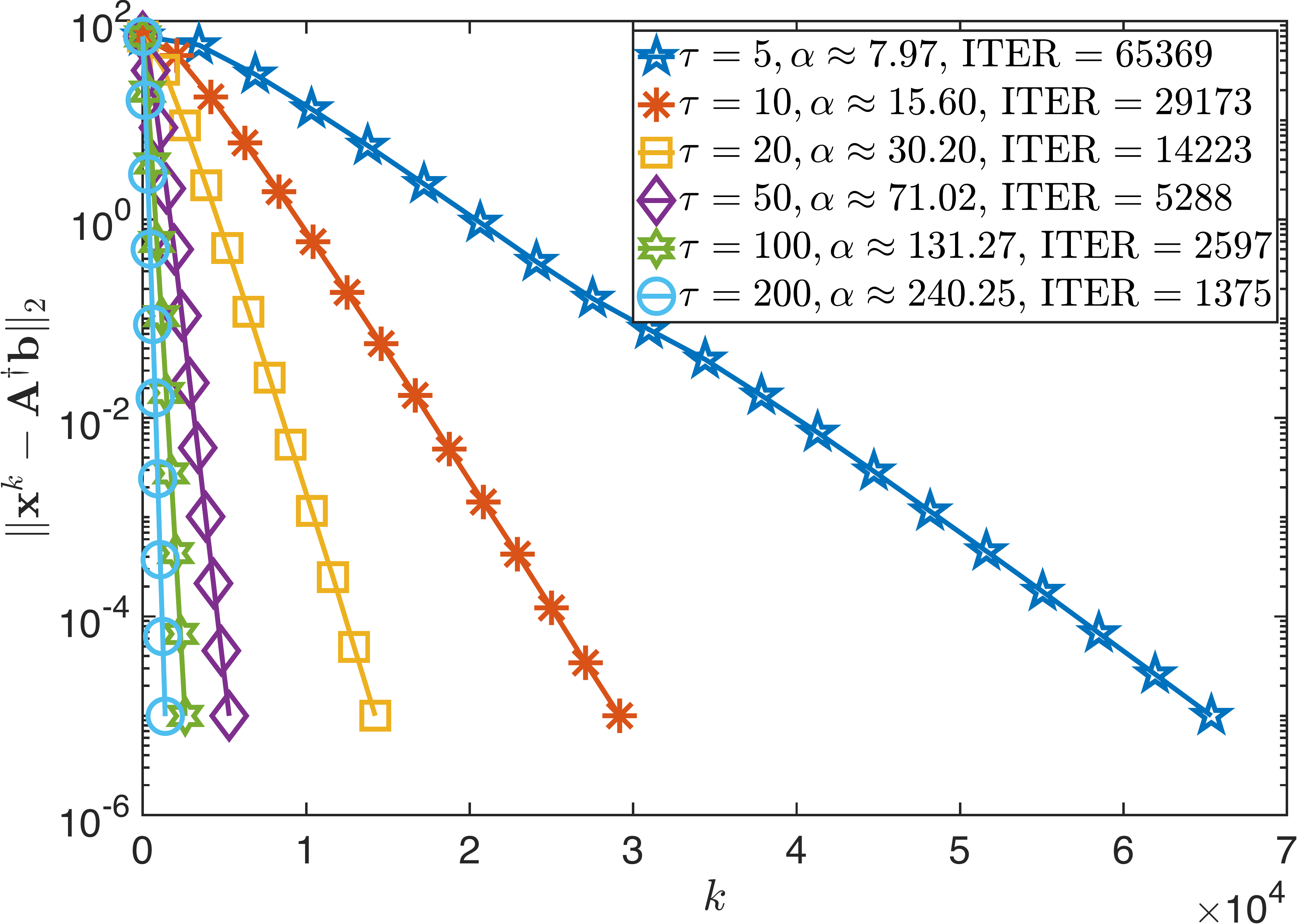,height=2.15in}\quad\epsfig{figure=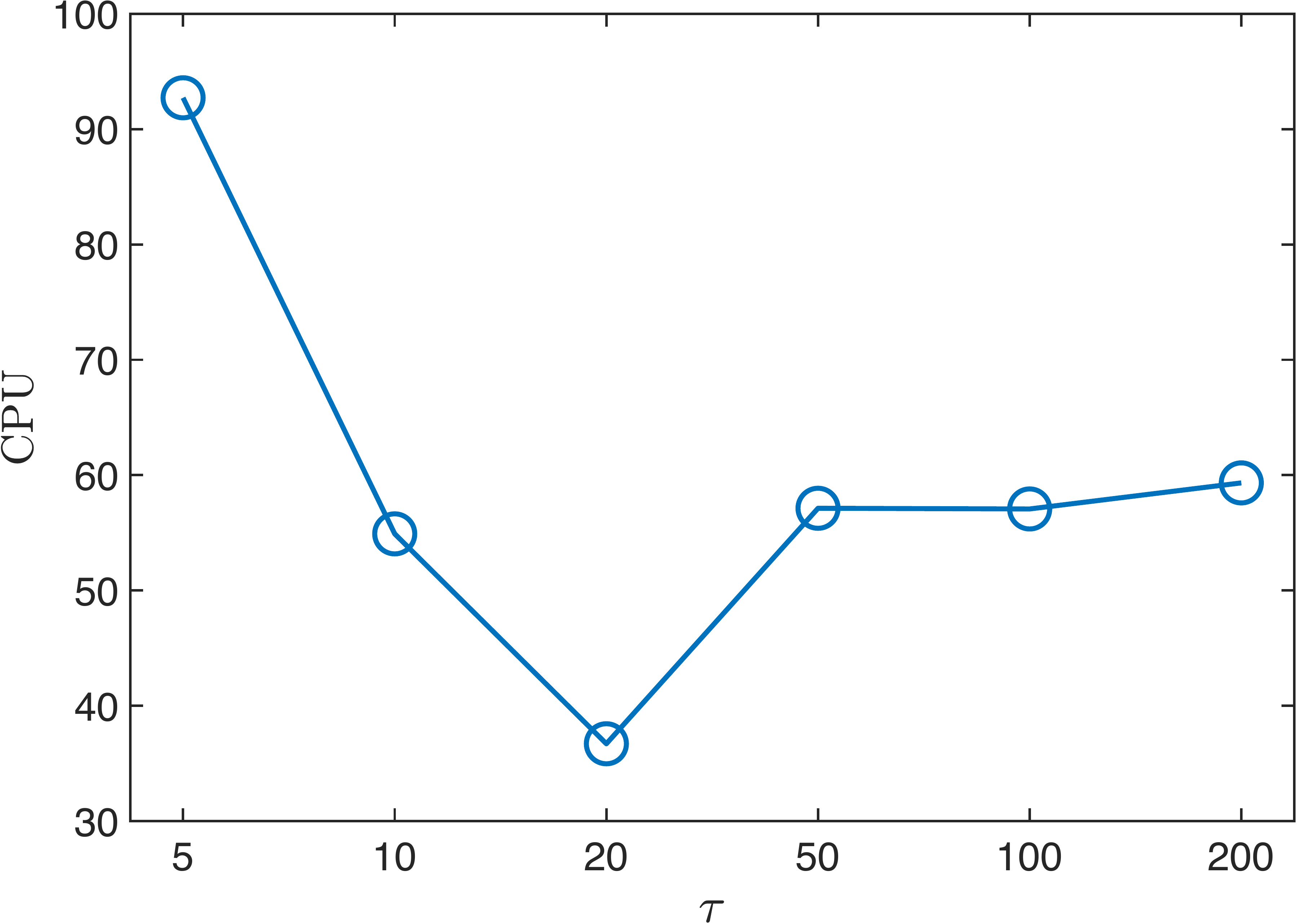,height=2.15in}}
\caption{The average (10 trials of each case) error $\|\mbf x^k-\mbf A^\dag\mbf b\|_2$ and CPU of REBK with different block sizes $\tau=10,50,100,200$ and stepsize $\alpha=1.75/\beta_{\max}$ for inconsistent linear systems. The average numbers of required iterations are also reported. Upper: Type I matrix $\bf A=UDV^\rmt$ with $m=20000$, $n=5000$, $r=4500$, and $\kappa=2$. Lower: Type II matrix $\bf A=$ {\tt randn(20000,5000)}.}\label{fig2}
\end{figure}

{ In Figure \ref{fig2}, we plot the error $\|\mbf x^k- \mbf A^\dag\mbf b\|_2$ and the computing times of REBK with block sizes $\tau=5,10,20,50,100,200$ and stepsize $\alpha=1.75/\beta_{\max}$ for two inconsistent linear systems  with coefficient matrices of Types I ($\bf A=UDV^\rmt$ with $m=20000$, $n=5000$, $r=4500$, $\kappa=2$) and II ($\bf A=$ {\tt randn(20000,5000)}). The average numbers of required  iterations are also reported. We observe: (i) increasing block size and using the empirical stepsize $\alpha=1.75/\beta_{\max}$ lead to a better convergence in terms of the numbers of iterations; (ii) with the increase of block size, the computing time first decreases, then increases after reaching the minimum value, and finally tends to be stable. This means that for sufficiently large block size the decrease in iteration complexity cannot compensate for the increase in cost per iteration. On the other hand, if a distributed version of REBK is implemented, a larger $\tau$ will be better.}

\begin{figure}[!htpb]
\centerline{\epsfig{figure=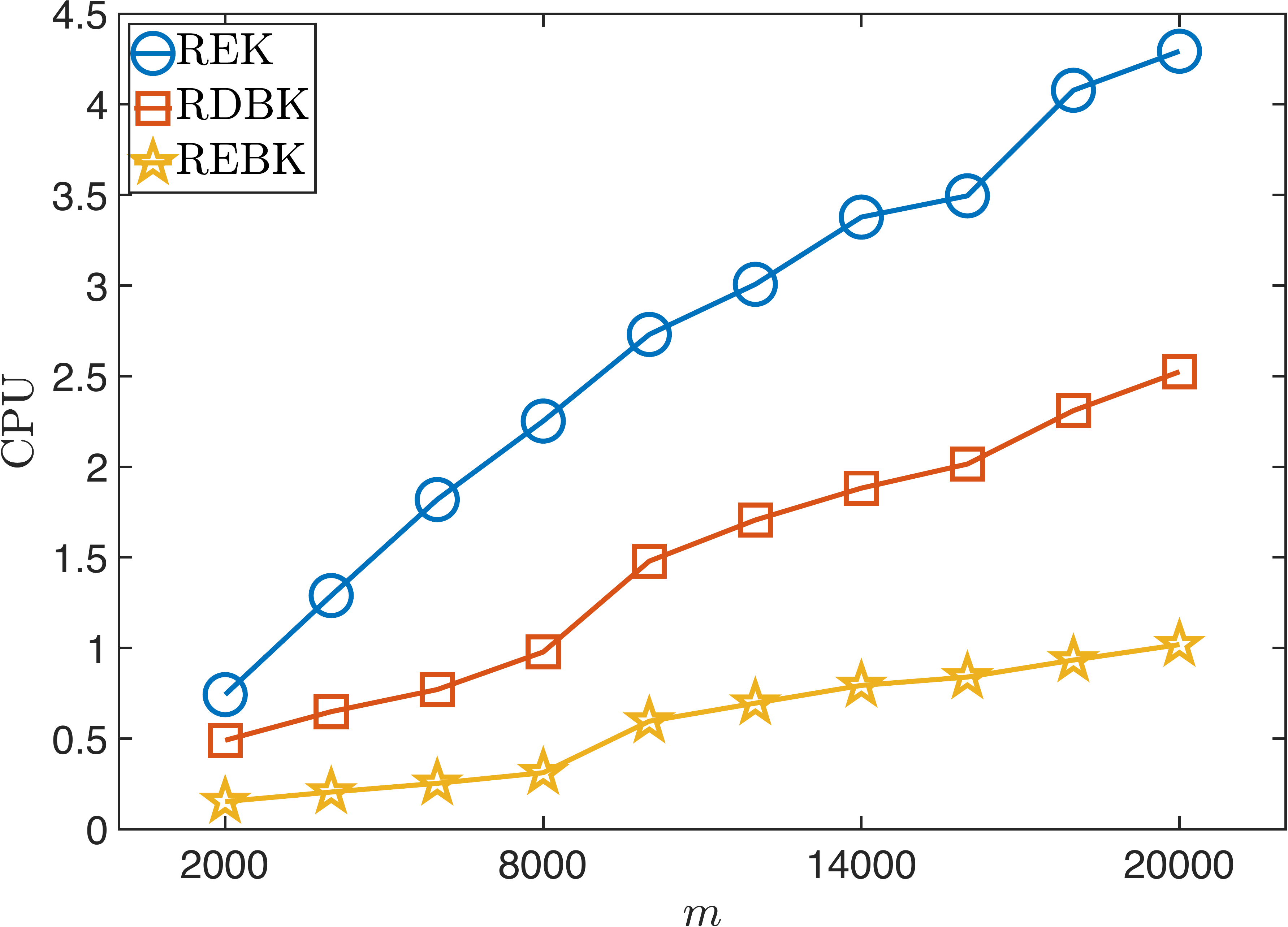,height=2.15in}\quad\epsfig{figure=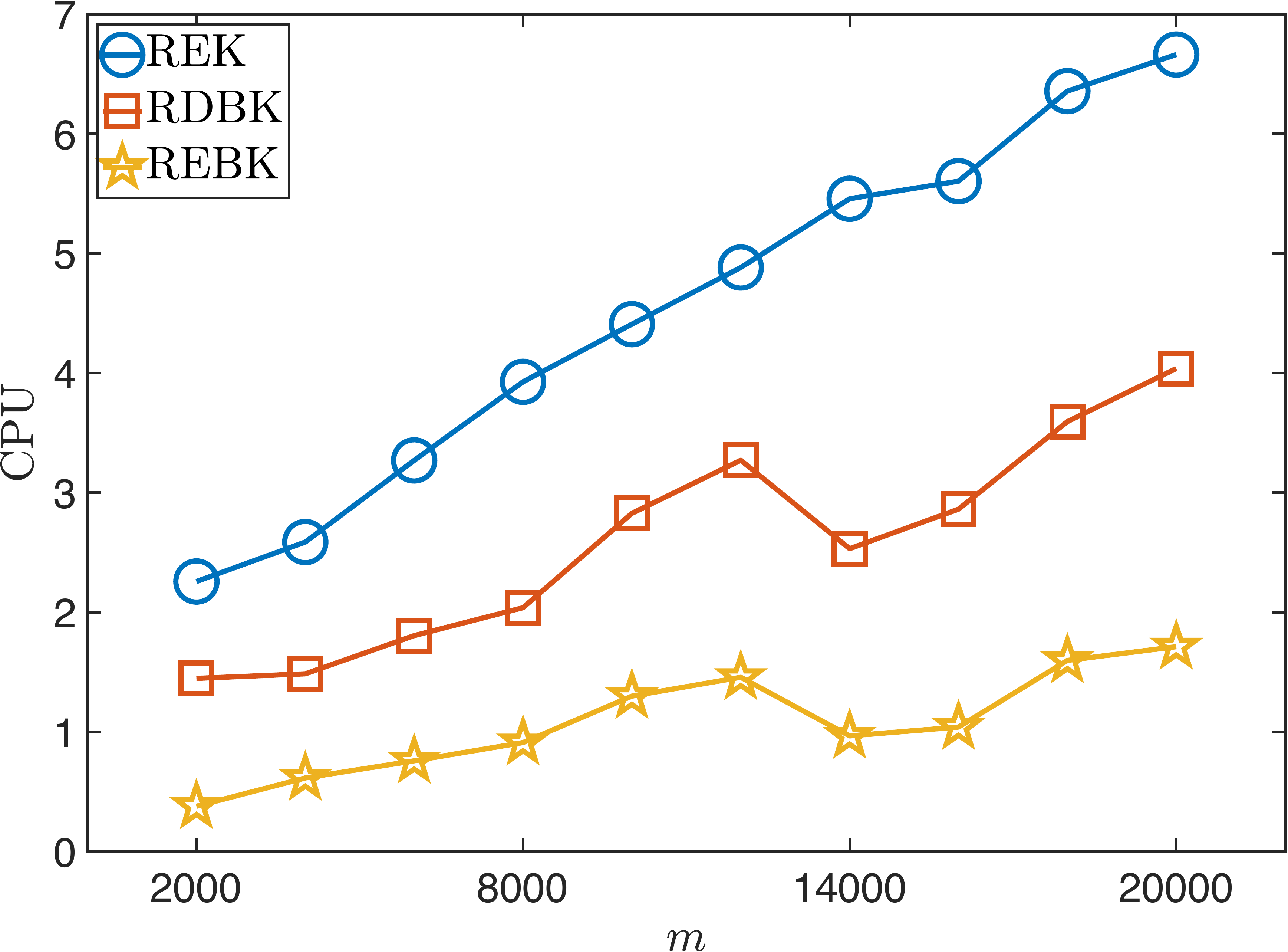,height=2.15in}}
\caption{The average (10 trials of each algorithm) CPU of REK, RDBK($\tau=10$), and REBK($\tau=10,\alpha=1.75/\beta_{\max}$) for inconsistent linear systems. Left: Type I matrix $\bf A=UDV^\rmt$ with $m=2000,4000,\ldots,20000$, $n=500$, $r=250$, and $\kappa=2$. Right: Type II matrix $\bf A=$ {\tt randn(m,n)} with $m=2000,4000,\ldots,20000$ and $n=500$.}\label{fig3}
\end{figure}

{In Figure \ref{fig3}, we plot the computing times of the REK, RDBK, and REBK algorithms for inconsistent linear systems with coefficient matrices of Types I ($\bf A=UDV^\rmt$ with $m=2000$, $4000$, $\ldots$, $20000$, $n=500$, $r=250$, $\kappa=2$) and II ($\bf A=$ {\tt randn(m,n)} with $m=2000$, $4000$, $\ldots$, $20000$, $n=500$). For all cases, the block size $\tau=10$ and the stepsize $\alpha=1.75/\beta_{\max}$ are used. We observe that both RDBK and REBK are better than REK, and REBK is the best.}

\subsection{Real-world data}
Finally, we test REK and REBK {using eight inconsistent linear systems with coefficient matrices} from the University of Florida sparse matrix collection \cite{davis2011unive}. The eight matrices are {\tt abtaha1}, {\tt flower\_5\_1}, {\tt football}, {\tt lp\_nug15}, {\tt relat6}, {\tt relat7}, {\tt Sandi\_authors}, and {\tt WorldCities}. In Table \ref{t3}, we report the numbers of iterations and the computing times for the REK and REBK algorithms. { For each matrix, we tested two stepsizes of REBK, the first is $1/\beta_{\max}$ and the second is empirical.} We observe that REBK based on good choices of block size and stepsize significantly outperforms REK. Moreover, good stepsize and block size are problem dependent.

\begin{table}[htp]
\caption{The average (10 trials of each algorithm) {\rm ITER} and {\rm CPU} of {\rm REK} and {\rm REBK}($\tau,\alpha$) for inconsistent linear systems with coefficient matrices from \cite{davis2011unive}. For each matrix, two stepsizes of {\rm REBK} are tested: the first is $1/\beta_{\max}$, and the second is empirical.}
\label{t3}
\begin{center} \footnotesize
\begin{tabular}{c|c|c|c|c|c|c|c|c|c|c} \toprule
\multirow{2}{*}{Matrix} & \multirow{2}{*}{$m\times n$}& \multirow{2}{*}{rank} & \multirow{2}{*}{$\dsp\frac{\sigma_1(\mbf A)}{\sigma_r(\mbf A)}$} & \multicolumn{2}{|c|}{REK} & \multicolumn{5}{|c}{REBK}\\
  \cline{5-11}&&&&ITER &CPU &$\tau$ & $\alpha$ & ITER &CPU &speed-up\\ \hline
\multirow{2}{*}{\tt abtaha1} & \multirow{2}{*}{$14596\times209$} & \multirow{2}{*}{209} & \multirow{2}{*}{12.23} & \multirow{2}{*}{276946} & \multirow{2}{*}{89.38} & \multirow{2}{*}{10} & 1.82 & 151395 &  68.30 & 1.31\\ \cline{8-11}&&&&&&&5 & 56064 &25.34 &3.53\\ \hline
\multirow{2}{*}{\tt flower\_5\_1}  & \multirow{2}{*}{$211\times201 $} & \multirow{2}{*}{179} & \multirow{2}{*}{13.70} & \multirow{2}{*}{135117} & \multirow{2}{*}{5.16} & \multirow{2}{*}{5} &1 & 136037  & 6.15 &  0.84 \\ \cline{8-11}&&&&&&& 4 &34381& 1.55 & 3.34 \\ \hline
\multirow{2}{*}{\tt football}  & \multirow{2}{*}{$35\times35 $} & \multirow{2}{*}{19} & \multirow{2}{*}{166.47}  & \multirow{2}{*}{810792} & \multirow{2}{*}{21.99} & \multirow{2}{*}{5} &1 & 858215 & 30.64 & 0.72 \\ \cline{8-11}&&&&&&& 2 & 409995 &14.63 &  1.50\\ \hline
\multirow{2}{*}{\tt lp\_nug15} & \multirow{2}{*}{$6330\times22275$} & \multirow{2}{*}{5698} &\multirow{2}{*}{2.73} & \multirow{2}{*}{216924} & \multirow{2}{*}{220.64} & \multirow{2}{*}{20} & 3.53 & 40539 & 199.67 & 1.10 \\  \cline{8-11}&&&&&&&5 &31039 & 158.29& 1.39\\ \hline
\multirow{2}{*}{\tt relat6}  & \multirow{2}{*}{$2340\times 157$} & \multirow{2}{*}{137} & \multirow{2}{*}{7.74} & \multirow{2}{*}{34536} & \multirow{2}{*}{2.43} & \multirow{2}{*}{10} & 1& 34273 & 3.81 & 0.64 \\\cline{8-11}&&&&&&&2.5 &13971 & 1.56 &1.56 \\ \hline
\multirow{2}{*}{\tt relat7}  & \multirow{2}{*}{$21924\times 1045$} & \multirow{2}{*}{1012} & \multirow{2}{*}{10.85} & \multirow{2}{*}{550810} & \multirow{2}{*}{283.69} & \multirow{2}{*}{10} & 1 & 542100 & 466.89 & 0.61 \\ \cline{8-11}&&&&&&&2.5 & 218287& 188.81 & 1.50\\ \hline
\multirow{2}{*}{\tt Sandi\_authors}  & \multirow{2}{*}{$86\times86 $} & \multirow{2}{*}{72}  & \multirow{2}{*}{189.58} & \multirow{2}{*}{2525141} & \multirow{2}{*}{73.28} & \multirow{2}{*}{5} & 1 & 2533343 & 99.36 & 0.74 \\ \cline{8-11}&&&&&&&2.5 & 999294& 39.15 & 1.87 \\ \hline
\multirow{2}{*}{\tt WorldCities}  & \multirow{2}{*}{$315\times 100$} & \multirow{2}{*}{100} & \multirow{2}{*}{66.00}& \multirow{2}{*}{120699} & \multirow{2}{*}{4.32} & \multirow{2}{*}{5} & 1.13 & 105647 & 4.52 & 0.96 \\ \cline{8-11}&&&&&&&2.5 & 47372 & 2.02& 2.14 \\  
   \bottomrule
\end{tabular}
\end{center}
\end{table}

\section{Concluding remarks}

We have proposed a  randomized extended block Kaczmarz (REBK) algorithm for solving general linear systems and prove its convergence theory. At each step, REBK uses two RABK (with special choice of weights) updates. The new algorithm can utilize efficient implementations on distributed computing units. Numerical experiments show that the crucial point for guaranteeing fast convergence is to obtain good block size and stepsize. Finding appropriate variable stepsize by the adaptive extrapolation \cite{necoara2019faste} and proposing more effective partitions based on the techniques of \cite{needell2014paved,drummond2015parti,torun2018ovel,necoara2019faste} should be valuable topics. { We also note that RABK allows the flexibility that the distributions from which blocks are selected do not require the blocks to form a partition of the columns, or rows. Designing variants of REBK based on RABK with random samplings that do not depend on the partitions is straightforward. We believe the technique used in the proof of Theorem \ref{main} still works for these variants. Although the analysis will be more complicated. Besides, developing parallel and accelerated variants of REBK based on the approach used by Richt\'{a}rik and Tak\'{a}\v{c} \cite{richtarik2017stoch} is also worth exploring. We will work on these topics in the future.}

\section*{Acknowledgments}
{The authors are thankful to the referees for their detailed comments and valuable suggestions that have led to remarkable improvements.} The research of the first author was supported by the National Natural Science Foundation of China (No.11771364) and the Fundamental Research Funds for the Central Universities (No.20720180008).

\end{document}